# CHAMBER STRUCTURE OF DOUBLE HURWITZ NUMBERS

RENZO CAVALIERI, PAUL JOHNSON, AND HANNAH MARKWIG

ABSTRACT. Double Hurwitz numbers count covers of $\mathbb{P}^1$ by genus $g$ curves with assigned ramification profiles over $0$ and $\infty$, and simple ramification over a fixed branch divisor. Goulden, Jackson and Vakil have shown double Hurwitz numbers are piecewise polynomial in the orders of ramification ([GJV05]), and Shadrin, Shapiro and Vainshtein have determined the chamber structure and wall crossing formulas for $g = 0$ ([SSV08]). This paper gives a unified approach to these results and strengthens them in several ways — the most important being the extension of the results of [SSV08] to arbitrary genus.

The main tool is the authors' previous work ([CJM]) expressing double Hurwitz number as a sum over certain labeled graphs. We identify the labels of the graphs with lattice points in the chambers of certain hyperplane arrangements, which are well known to give rise to piecewise polynomial functions. Our understanding of the wall crossing for these functions builds on the work of Varchenko ([Var87]), and could have broader applications.

## Contents

The second author is supported in part by NSF grants DMS-0602191 and DMS-0902754. The third author is supported by the German Research Foundation (Deutsche Forschungsgemeinschaft (DFG)) through the Institutional Strategy of the University of Göttingen and by the DFG grant MA 4797/1-1.

# 1. Introduction

## 1.1. Statement of Results.
Hurwitz theory studies holomorphic maps between Riemann surfaces with specified ramification. Double Hurwitz numbers count covers of $\mathbb{P}^1$ with assigned ramification profiles over $0$ and $\infty$, and simple ramification over a fixed branch divisor.

We use a new notation for double Hurwitz numbers. We define $H_g(\mathbf{x})$ to be the genus $g$ double Hurwitz number with profile $\mathbf{x_0} := \{x_i | x_i > 0\}$ over zero and $\mathbf{x}_\infty := \{x_i | x_i < 0\}$ over $\infty$ (previous notation recorded the ramification over $0$ and $\infty$ separately). Furthermore, frequently the natural numerical invariant is $r$, the number of simple ramifications, rather than the genus $g$. Since these are equivalent by the Riemann-Hurwitz formula, we use $H^r(\mathbf{x})$ to denote $H_g(\mathbf{x})$ when it makes formulas more attractive.

Our first result is a new proof of the following theorem in [GJV05]:

**Theorem 1.1** (GJV). *The function*
$$H_g(\mathbf{x}) : \big\{ \sum x_i = 0 \big\} \subset (\mathbb{Z} \setminus \{0\})^n \to \mathbb{Q}$$
*is a piecewise polynomial function of degree $4g - 3 + n$.*

Our techniques allow us to answer an (implicit - see Section 1.4) conjecture of Goulden, Jackson and Vakil.

**Theorem 1.2.** *$H_g(\mathbf{x})$ is either even or odd, depending on the parity of the leading degree $4g - 3 + n$.*

We extend the results of [SSV08] to all genera. First, we determine the regions on which $H_g(\mathbf{x})$ is polynomial:

**Theorem 1.3.** *The chambers of polynomiality of $H_g(\mathbf{x})$ are bounded by walls corresponding to the resonance hyperplanes $W_I$, given by the equation*
$$W_I = \left\{ \sum_{i \in I} x_i = 0 \right\},$$
*for any $I \subset \{1, \ldots, n\}$.*

Our main result is a wall crossing formula. We denote the chambers of the resonance arrangement as $H$-*chambers*.

**Definition 1.4.** Let $\mathfrak{c}_1$ and $\mathfrak{c}_2$ be two $H$-chambers adjacent along the wall $W_I$, with $\mathfrak{c}_1$ being the chamber with $x_I < 0$. The Hurwitz number $H^r(\mathbf{x})$ is given by polynomials, say $P_1(\mathbf{x})$ and $P_2(\mathbf{x})$, on these two regions. A wall crossing formula is a formula for the polynomial
$$WC_I^r(\mathbf{x}) = P_2(\mathbf{x}) - P_1(\mathbf{x}).$$

Note that with the notation $WC_I^r(\mathbf{x})$ there is no ambiguity about which direction we cross the wall.

**Theorem 1.5.**

$$
(1) \quad WC_I^r(\mathbf{x}) = \sum_{\substack{s+t+u=r \\ |\mathbf{y}|=|\mathbf{z}|=|\mathbf{x}_I|}} \left( (-1)^t \cdot \binom{r}{s,t,u} \cdot \frac{\prod \mathbf{y}_i}{\ell(\mathbf{y})!} \cdot \frac{\prod \mathbf{z}_j}{\ell(\mathbf{z})!} \cdot \right.
$$
$$
\left. H^s(\mathbf{x}_I, \mathbf{y}) \cdot H^{t\bullet}(-\mathbf{y}, \mathbf{z}) \cdot H^u(\mathbf{x}_{I^c}, -\mathbf{z}) \right)
$$

Here $\mathbf{y}$ is an ordered tuple of $\ell(\mathbf{y})$ positive integers with sum $|\mathbf{y}|$, and similarly with $\mathbf{z}$.

This formula appears not to depend on the particular choice of chambers $\mathfrak{c}_1$ and $\mathfrak{c}_2$ that border on the wall, but only upon the wall $W_I$; however the polynomials for the simpler Hurwitz numbers in the formula depend on chambers themselves.

The walls $W_I$ correspond to those values of $\mathbf{x}$ where the cover could potentially be disconnected, or where $\mathbf{x}_i = 0$. Crossing this second type of wall corresponds to moving a ramification between 0 and $\infty$. In the traditional view of double Hurwitz numbers, this would cross between different problems: the lenght of the profiles over 0 and $\infty$ were fixed separately, rather than just the sum of the lengths being fixed. However, Theorem 1.5 suggests that it is natural to treat these as part of the same problem: the wall crossing formula for $\mathbf{x}_i = 0$ is identical to the other wall crossing formulas.

1.2. **Overview of Methods.** This paper is an exploration of the consequences of formula (2) in the authors' previous work [CJM], which expresses double Hurwitz numbers $H_g(\mathbf{x})$ as a sum over certain graphs $\Gamma$, which we call monodromy graphs (Definition 2.1). Each internal edge of a monodromy graph is labeled with a positive integer. The contribution of each monodromy graph is the product of these integers. In genus zero, these edge labelings are determined uniquely by $\mathbf{x}$, and as a result the genus zero case of all of our theorems follow quickly from the graphs, as is presented in Section 6 of [CJM].

In positive genus, however, the choice of edge labelings is not unique. The crux of this paper is to understand the space of edge labelings for each directed graph and value of $\mathbf{x}$. We show in Section 2 that for each directed graph and value of $\mathbf{x}$, the space of edge labelings (which we call *flows*) are the lattice points in certain bounded polytopes we call $F$-chambers. As $\mathbf{x}$ changes, the faces of the $F$-chambers shift, but their normal directions remain constant. Thus, for each directed graph the contribution is the sum of a polynomial (the product of the edge weights) over the lattice points in a polytope. Furthermore, the $F$-chambers for directed graphs with the same underlying undirected graph $\Gamma$ fit together as the set $\mathcal{BC}_\Gamma(x)$ of bounded chambers of a natural hyperplane arrangement $\mathcal{A}_\Gamma(\mathbf{x})$ associated to $\Gamma$ and $\mathbf{x}$. Intuitively, the arrangements $\mathcal{A}_\Gamma(\mathbf{x})$ have an easy description: we consider flows of water along the edges of the graph. At vertices not labeled by a part of $\mathbf{x}$, water is conserved, and the inflowing water must equal the outflowing water. The vertices labeled by parts of $\mathbf{x}$, however, have valves, that add in $\mathbf{x}_i$ of water. The space of all such flows is the underlying affine space of the arrangement, and the hyperplanes are given by when the flow along a given edge is equal to zero.

There is a general theory of lattice points in polytopes that can be brought to bear upon the problem. Theorems 1.1, 1.2 and 1.3 follow from standard results in this theory. Since for all integral $\mathbf{x}$ the vertices of the flow polytopes are integers, we can conclude that the sums are piecewise polynomial, and that the walls occur when the topology of the hyperplane arrangement changes, proving Theorems 1.1 and 1.3.

Theorem 1.2 is more subtle: as opposed to integrating a homogenous polynomial over a polytope, summing a homogenous polynomial over the lattice points of a polytope does not in general result in an odd or even polynomial. However, Ehrhart reciprocity says that the failure to be odd or even is essentially due to the lattice points contained in the boundary of the polytope. In our case the polynomial we are summing vanishes on the boundary of the polytope, giving 1.2. This phenomenon plays an important role in our approach to Theroem 1.5.

1.3. **Wall Crossing.** The wall crossing phenomenon for polytopes is a rich area with many possible approaches and ongoing research. In our proof of Theorem 1.5, we follow Varchenko's viewpoint in [Var87]; however, other approaches should also prove fruitful in investigating this problem – forthcoming work of Ardila [Ard] takes the point of view of generalized Dahmen-Micchelli spaces [dPV]. As Varchenko's viewpoint is less standard than the other lattice point techniques we use, we give a short overview.

The fundamental idea is that understanding wall crossing becomes much simpler if instead of focusing on a single polytope, we view the polytope as one chamber of a hyperplane arrangement. Then wall crossing for a given polytope can be understood in terms of adding and subtracting certain nearby polytopes in the arrangement. The data of which polytopes to add and subtract is encoded in a linear map called the Gauss-Manin connection. For the explanation of the name, see

Section 5. We now illustrate this idea in the local setting. A global example is worked out in detail in Example 4.1.

Consider the 1-dimensional family $\mathcal{A}_t$ of $n$-dimensional arrangements, where the $n+1$ hyperplanes of $\mathcal{A}_t$ consist of the $n$ coordinate hyperplanes, together with the hyperplane $x_1 + \cdots + x_n = t$. For $t > 0$, the arrangement is simple, and there are $2^{n+1} - 1$ chambers: one bounded $n$-simplex and $2^{n+1} - 2$ unbounded chambers, each of which borders the bounded chamber on one of its proper faces. As $t$ approaches zero, the bounded chamber shrinks, until it disappears at $t = 0$ and there is a nontransverse intersection: all $k+1$ hyperplanes intersect at the origin. When $t < 0$, a new bounded simplex $A$ appears. Furthermore, the topology of each of the unbounded chambers has changed, but in an easily described way: we must add or subtract the appearing chamber $A$ to each unbounded chamber $U$, depending on the codimension with which $A$ and $U$ border.

If we add a few fixed hyperplanes to this hyperplane arrangement to obtained bounded chambers, the volume of one of the resulting chambers would be a polynomial in $t$ for $t > 0$ and for $t < 0$; the difference between these polynomials is, up to a sign, the volume of the appearing/vanishing simplex. This is essentially the local picture for all wall crossings in families where the generic arrangement is simple: at the wall, there are $k+1$ hyperplanes meeting in codimension $k$. On either side of the wall, these $k+1$ hyperplanes bound a $k$-simplex crossed with $\mathbb{R}^{n-k}$, which may be further cut into smaller chambers by the other hyperplanes. We call these chambers vanishing chambers. As we cross the wall, to keep using the same volume polynomial, we must add or a subtract each vanishing/appearing chamber to each of the $2^{k+1} - 2$ chambers obtained by crossing some proper subset of the $k+1$ hyperplanes that meet nontransversely.

Varchenko's approach to understanding this phenomenon is to use cones: any chamber can be written as a signed sum of cones, and cones obviously do not change topology - rather, the change in topology is due to the change in relative position of the cones.

Care is needed in extending this approach to integer points in the polytopes: for each boundary face, we must specify whether we are including the lattice points on that face or not. This can be done from the cone point of view: when writing our polytope as a sum of cones, we must specify whether each face of the cone is included or not. If we keep track of this information, everything follows. Varchenko calls the result a "combinatorial connection": it is the usual Gauss-Manin connection, with corrections by lower dimensional cells. Thus, the combinatorial connection can be understood as a generalization of Ehrhart reciprocity([BR07]).

We wish to apply the general machinery of the combinatorial connection to our situation. As the polynomial we are summing over the lattice points vanishes on the boundary, the lower dimensional corrections to the combinatorial connection are unnecessary for our purposes, and we only need a formula for the usual connection. The difficulty is that the generic hyperplane arrangement in our families is not simple, and so the easy description of the Gauss-Manin connection described above fails.

The technical heart of our paper is a formula for the Gauss-Manin connection for the families of hyperplane arrangements we are dealing with, stated in terms of the combinatorics of the graph (Lemma 6.12).

1.4. **Connections with geometry.** Although our methods are essentially combinatorial, much of the motivation for studying Hurwitz theory lies in deep geometric connections to moduli spaces. We now discuss the relationship of our work to these results.

Historically, one motivation for studying Hurwitz theory is to understand the moduli space of curves $\mathcal{M}_g$ - in particular, to show it is irreducible. Another recent connection is ELSV formula [ELSV01], which expresses single Hurwitz numbers – i.e., when there is no ramification over $\infty$ – in terms of intersection numbers on $\overline{\mathcal{M}}_{g,n}$. The role of the ELSV formula in understanding these and related intersection numbers has been remarkably fruitful, and a survey would be far beyond our needs. We mention one result going in the opposite direction: the geometric form of the ELSV formula proves and explains polynomiality for single Hurwitz numbers as conjectured by Goulden, Jackson and Vainshtein [GJV00].

Goulden, Jackson and Vakil conjecture in [GJV05] that there should be a formula similar to the ELSV formula for one part double Hurwitz numbers – those Hurwitz numbers with total ramification over 0 and arbitrary ramification over $\infty$. One part double Hurwitz numbers are in fact polynomial, and the GJV conjecture would provide a geometric explanation why, parallel to that for single Hurwitz numbers. In their conjecture, the moduli space of curves is replaced by some universal Picard scheme which over the smooth locus parameterizes a complex curve together with a line bundle; the difficulty is determining how to compactify the Picard group of nodal curves.

One point where our work makes interesting contact with this conjecture is our observation that double Hurwitz numbers behave well when we allow ramification to become negative and pass from zero to infinity. From this perspective, there is nothing special about one part double Hurwitz numbers. They are one chamber of polynomiality for $H_g(\mathbf{x})$ – where $\mathbf{x_1} > 0$ and $\mathbf{x_i} < 0$ for $i > 1$.

It is then natural to wonder whether the GJV conjecture could be extended to give a geometric explanation for the polynomiality of each chamber, and indeed the whole piecewise polynomial nature of $H_g(\mathbf{x})$. Ideally, one would hope that there were stability conditions that lead to different compactifications of the Picard varieties of nodal curves, and that the changing choice in compactification would account for the change in polynomial.

Part of the motivation for the GJV conjecture was a rather explicit formula for one part double Hurwitz numbers obtained via representation theory; this approach is extended to all chambers in [Joh], providing another proof of Theorems 1.1, 1.2 and 1.3, as well as a proof of the strong piecewise polynomiality conjecture of [GJV05] that the methods of the current paper cannot prove (see Section 3.4). Taken together, these results show that the algebraic form of all double Hurwitz polynomials are compatible with an extended GJV conjecture.

Though all three approaches to Theorem 1.1 are largely algebraic, it is interesting to observe that both the original proof in [GJV05] and the proof presented here have elements that point toward connections with the moduli space of curves $\mathcal{M}_{g,n}$: the proof in [GJV05] uses ribbon graphs, which index the cells of a combinatorial description of $\mathcal{M}_{g,n}$ (for an introduction, see [LZ04]). The trivalent graphs that we use index the top dimensional cells of the stacky fan that is the closest thing to a tropical $\overline{\mathcal{M}}_{g,n}$ (see [BMV]). This is not surprising, as our original motivation was from tropical geometry.

Our method provides some evidence for an extended GJV conjecture as described above: the hyperplane arrangements $\mathcal{A}_\Gamma(\mathbf{x})$ appearing on our work are precisely the combinatorial information used by Oda and Seshadri [OS79] (see also [Ale04]) to construct their compactifed Jacobians $Jac_\phi$ for nodal curves. In their work, $\Gamma$ is the dual graph of a nodal curve, they are working considering the infinite hyperplane arrangement given by all integer translates of our hyperplanes, and the polyhedra are used to construct toric varieties; in any case, there is a chamber structure on possible values of $\phi$ given by the changing topology of the arrangement. Thus, the combinatorial mechanism responsible for producing our chamber structure is identical to that producing changing stability conditions for compactified Jacobians.

Nice behavior as ramification crosses from 0 to $\infty$ has been observed previously in the closely related area of the Gromov-Witten theory of $\mathbb{P}^1$ relative to 0 and $\infty$. In [OP06], Okounkov and Pandharipande show that the one point invariants are polynomial of degree $2g$ in the orders of ramification, and observe, in what they call "crossing symmetry", that these polynomials are completely symmetric. Additionally, this phenomenon has been observed by Vakil ([Vak]) for the pushforward of the rubber virtual fundamental class to Chow classes in the moduli space $M_g^{rt}$ of curves with rational tails.

| | | |
|---|---|---|
| $\mathbf{x}$ | ordered tuple of $n$ nonzero integers $x_i$ with $\sum_i x_i = 0$ | p. 2 |
| $H_g(\mathbf{x}), H^r(\mathbf{x})$ | Hurwitz numbers ($g$ genus, $r$ #simple ramifications) | p. 2 |
| $W_I$ | (Hurwitz) walls | 1.3 |
| $WC_I^r(\mathbf{x})$ | wall crossing | 1.4 |
| $\mathfrak{c}$ | A (Hurwitz) chamber for $\mathbf{x}$ | 1.4, 3.5 |
| $\Gamma(d,o), \Gamma(\mathbf{x},d,o)$ | monodromy graph (directed $\mathbf{x}$-graph with a compatible ordering of vertices) | 2.1 |
| $\varphi_\Gamma$ | product of internal weights of $\Gamma$ | (2) |
| $\Gamma, \Gamma(\mathbf{x})$ | $\mathbf{x}$-graph | 2.3 |
| $\Gamma(d), \Gamma(\mathbf{x},d)$ | directed $\mathbf{x}$-graph | 2.3 |
| $F_\Gamma(\mathbf{x})$ | space of flows on $\Gamma(\mathbf{x})$ | 2.8 |
| $\mathcal{A}, \mathcal{A}_\Gamma(\mathbf{x})$ | flow hyperplane arrangement | 2.8 |
| $\varphi_\mathcal{A}$ | defining equation of $\mathcal{A}$ | 2.8 |
| $A, B, \dots$ | $F$-chambers | 2.15 |
| $\Gamma_A, \Gamma_B, \dots$ | directed $\mathbf{x}$-graphs associated to the $F$-chambers $A, B, \dots$ | 2.15 |
| $m(\Gamma_A), m(A)$ | number of vertex orderings of $\Gamma_A$ | 2.15 |
| $\text{sign}(A)$ | $(-1)^e$ where $e = $ #edges in which $\Gamma_A$ differs from reference orientation | 2.15 |
| $\text{Ch}(\mathcal{A}_\Gamma(\mathbf{x}))$ | set of chambers of $\mathcal{A}_\Gamma(\mathbf{x})$ | 2.15 |
| $\mathcal{D}$ | discriminant arrangement | 3.4 |
| $\mathcal{BC}_\Gamma(\mathbf{x})$ | set of bounded $F$-chambers of $\mathcal{A}_\Gamma(\mathbf{x})$ | p. 32 |
| $\nabla_{\Gamma,12}, \nabla^*_{\Gamma,12}$ | Gauss-Manin connection resp. adjoint | 5.1 |
| $\Gamma/E$ | contracting the edges in $E$ | 6.3 |
| $C$ | a cut | 6.4 |
| $\text{Cuts}_I(\Gamma)$ | poset of cuts | 6.4 |
| $\Gamma'_A$ | contracting everything except the maximal cut | 6.7 |
| $X_{\text{Cuts}_I(\Gamma)}$ | cone of cuts | 6.7 |
| $X_C$ | face of $X_{\text{Cuts}_I(\Gamma)}$ corresponding to a cut | 6.7 |
| $L(P)$ | face lattice of the polyhedron $P$ | 6.8 |
| $\mathcal{K}$ | a cone | p. 42 |
| $\mathcal{K}(\mathfrak{c})$ | sum of chambers in a cone | p. 42 |
| $\mathcal{P}$ | partial orientation of edges | p. 42 |
| $\mathcal{K}_\mathcal{P}$ | cone corresponding to $\mathcal{P}$ | p. 43 |
| $\mathcal{O}_\Gamma$ | set of all orientations of $\Gamma$ | 7.2 |
| $\nabla^\mathcal{O}_{\Gamma,12}$ | graph connection | 7.2 |
| $\mathcal{NG}_\Gamma(\mathfrak{c})$ | nongeometric orientations | 7.2 |
| $\mathcal{K}^\mathcal{O}_\mathcal{P}$ | a combinatorial cone | 7.2 |
| $t(C)$ | thin cut associated to $C$ | 8.1 |
| $\gamma(T)$ | the middle components | 8.3 |

TABLE 1. Notation used throughout the paper.

1.5. **Organization.** The organization of this paper is as follows. Section 2 begins with a motivating example before introducing our interpretation of double Hurwitz numbers in terms of the hyperplane arrangements $\mathcal{A}_\Gamma$. We then apply this interpretation in Section 3 to prove Theorems 1.1, 1.2 and 1.3.

The rest of the paper is devoted to proving our wall crossing formula, Theorem 1.5. Section 4 is a gentle, example-oriented introduction and overview of our approach.

Section 5 formally discusses the combinatorial Gauss-Manin connection. In Section 6, we introduce the poset of cuts, which allows us to state our formula for the Gauss-Manin connection, and show how this formula implies the "heavy" wall crossing formula. Section 7 then proves our formula for the Gauss-Manin connection. In Section 8 we show how the main wall crossing formula follows from the "heavy" one.

1.6. **Acknowledgements.** We would like to thank Federico Ardila and Michael Shapiro for many helpful discussions.

## 2. Monodromy Graphs and Hyperplane Arrangements

Our main tool is the key observation of our earlier paper [CJM], that the cut and join recursion can conveniently be organized in terms of certain graphs; we review this in 2.1. We refine this organization by introducing certain hyperplane arrangements $\mathcal{A}_\Gamma(\mathbf{x})$. This is motivated with examples in 2.2 and presented formally in 2.3.

2.1. **Hurwitz numbers and Monodromy graphs.** The double Hurwitz number $H_g(x_1, \ldots, x_n)$ counts the number of maps $\pi : C \to \mathbb{P}^1$, where $C$ is a connected, genus $g$ curve and $\pi$ has profiles $\mathbf{x_0} := \{x_i | x_i > 0\}$ (resp. $\mathbf{x_\infty} := \{x_i | x_i < 0\}$) over $0$ (resp. $\infty$), and simple ramification over $r = 2g - 2 + n$ fixed other points. The preimages of $0$ and $\infty$ are marked. Furthermore, each cover is counted with weight $1/|\mathrm{Aut}(\pi)|$.

In [CJM, Lemma 4.1], we associate to each cover $\pi$ as above a decorated graph $\Gamma(=\Gamma(\mathbf{x}, d, o))$ that we call a monodromy graph.

**Definition 2.1.** For fixed $g$ and $\mathbf{x} = (x_1, \ldots, x_n)$, a graph $\Gamma$ is a *monodromy graph* if:

(1) $\Gamma$ is a connected, genus $g$, directed graph.
(2) $\Gamma$ has $n$ 1-valent vertices called *leaves*; the edges leading to them are *ends*. All ends are directed inward, and are labeled by the weights $x_1, \ldots, x_n$. If $x_i > 0$, we say it is an *in-end*, otherwise it is an *out-end*.
(3) All other vertices of $\Gamma$ are 3-valent, and are called *internal vertices*. Edges that are not ends are called *internal edges*.
(4) After reversing the orientation of the out-ends, $\Gamma$ does not have sinks or sources[1].
(5) The internal vertices are ordered compatibly with the partial ordering induced by the directions of the edges.
(6) Every internal edge $e$ of the graph is equipped with a *weight* $w(e) \in \mathbb{N}$. The weights satisfy the *balancing condition* at each internal vertex: the sum of all weights of incoming edges equals the sum of the weights of all outgoing edges.

The notation $\Gamma(\mathbf{x}, d, o)$ indicates that the graph comes with directed edges ($d$) and with a compatible vertex ordering ($o$).

*Remark* 2.2. Since the vertices in a monodromy graph are totally ordered, any orientation occurring has no directed cycles.

It follows from [CJM, Lemma 4.1] that the Hurwitz number is computed as:

$$(2) \qquad H_g(\mathbf{x}) = \sum_{\Gamma} \frac{1}{|Aut(\Gamma)|} \varphi_\Gamma,$$

where the sum is over all monodromy graphs $\Gamma$ for $g$ and $\mathbf{x}$, and $\varphi_\Gamma$ denotes the product of weights of all internal edges.

We simplify the combinatorics of this sum by grouping together families of monodromy graphs that coincide after forgetting structure.

**Definition 2.3.** Given $g$ and $\mathbf{x}$, an **x**-*graph* $\Gamma(\mathbf{x})$ (or simply $\Gamma$ if there is no risk of confusion) is a connected, genus $g$, trivalent graph with $n$ ends labeled $x_1, \ldots, x_n$.

*Remark* 2.4. Unless otherwise specified an **x**-graph is not a directed graph. In order to compare different directed graphs that map to the same **x**-graph it is useful to pick once and for all a reference orientation for all the edges. When assigning weights to the edges we understand that a positive weight preserves the reference orientation, whereas a negative weight reverses it. The convention that all ends in the reference orientation are directed inwards is compatible with positive ends being inputs and negative ends outputs.

---

[1] We do not consider leaves to be sinks or sources.

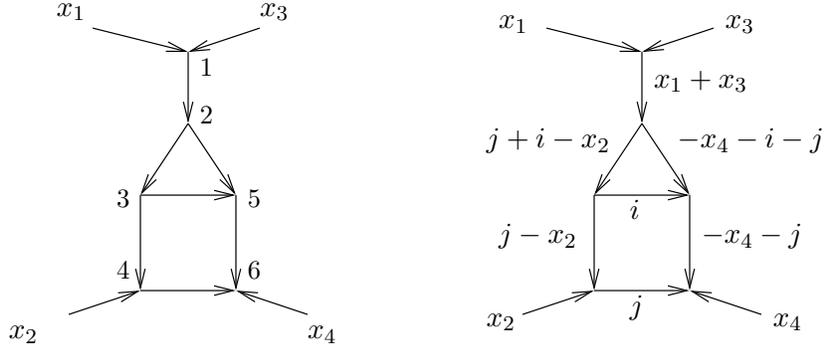

FIGURE 1. A directed **x**-graph and the weights of internal edges determined by the balancing condition.

The weights of internal edges for a directed **x**-graph are parameterized by (the integer points of) a $g$-dimensional polytope. The polytopes corresponding to different orientations of the edges of the same **x**-graph fit together as the bounded chambers of a natural hyperplane arrangement. In the remainder of this section we develop this point of view.

2.2. **A motivating example.**

**Example 2.5.** Consider the directed **x**-graph $\Gamma(\mathbf{x}, d, o)$ on the left hand side in Figure 1. We describe all monodromy graphs that equal $\Gamma(\mathbf{x}, d, o)$ after forgetting the weights of the internal edges. There are no monodromy graphs that equal $\Gamma(\mathbf{x}, d, o)$ after forgetting the weights if $x_1 + x_3 \leq 0$, so we assume that $x_1 + x_3 > 0$. Imposing the balancing condition at interior vertices leaves two degrees of freedom for the weights of interior edges, one for each independent cycle of $\Gamma$, as shown in the right hand side of Figure 1. All possible collections of edge labels are indexed by the lattice points in the polytope defined by:

$$\begin{aligned} i \geq 0, && j \geq 0, \\ j + i - x_2 \geq 0, && -x_4 - i - j \geq 0, \\ -x_4 - j \geq 0, && j - x_2 \geq 0, \end{aligned}$$

Figure 2 shows all hyperplanes $w(e) = 0$ with a normal vector indicating on which side of the hyperplane the inequality $w(e) > 0$ is satisfied.

The contribution of $\Gamma(\mathbf{x}, d, o)$ to $H_g(\mathbf{x})$ is given by

$$(x_1 + x_3) \cdot \sum i \cdot j \cdot (j + i - x_2) \cdot (-x_4 - i - j) \cdot (-x_4 - j) \cdot (j - x_2)$$

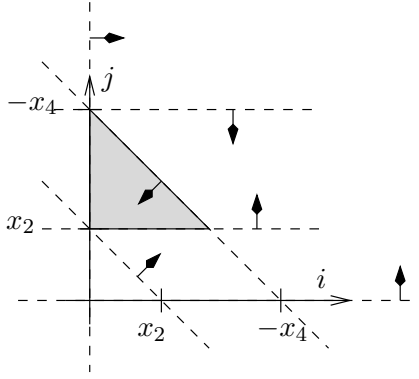

FIGURE 2. The polygon parameterizing internal edge weights for $\Gamma(\mathbf{x}, d, o)$.

where the sum goes over all lattice points $(i, j)$ in the polygon above (note that $\Gamma(\mathbf{x}, d, o)$ has no automorphisms). This equals

$(x_1 + x_3) \cdot$
$$\sum_{i=0}^{-x_4-x_2} \sum_{j=x_2}^{-i-x_4} i \cdot j \cdot (j + i - x_2) \cdot (-x_4 - i - j) \cdot (-x_4 - j) \cdot (j - x_2)$$

Expanding the sum, we observe it is an odd polynomial in the entries of $\mathbf{x}$ of degree $9 = 4g + n - 3$.

Let us point out the features of Example 2.5: for each cycle in the graph, there is one degree of freedom in choosing the labelings of the interior edges; each directed $\mathbf{x}$-graph together with a vertex ordering gives rise to a $g$ dimensional polytope whose integer points parametrize monodromy graphs; varying the vector $\mathbf{x}$ results in parallel translating the faces of the polytope. As long as the topology of the polytope remains the same, the contribution of a given directed $\mathbf{x}$-graph to the Hurwitz number is a polynomial of degree $4g - 3 + n$.

This polynomial does not depend on the vertex ordering. Thus, we can forget the vertex ordering and weight each directed $\mathbf{x}$-graph by an appropriate multiplicity $m$ (Definition 2.15).

**Example 2.6.** For the graph $\Gamma(\mathbf{x}, d)$, obtained by forgetting the vertex ordering in Example 2.5, the multiplicity $m(\Gamma(\mathbf{x}, d))$ equals 2: the vertices 4 and 5 can change their role.

If we forget more structure and just consider the underlying $\mathbf{x}$-graph $\Gamma(x)$, then the polytopes for the various choices of directions fit together nicely.

**Example 2.7.** Consider the **x**-graph $\Gamma(\mathbf{x})$ underlying Example 2.5. Retain the orientation of the edges in Figure 1 as a reference orientation, and the labels $w(e)$ for the internal edges obtained from the balancing condition:

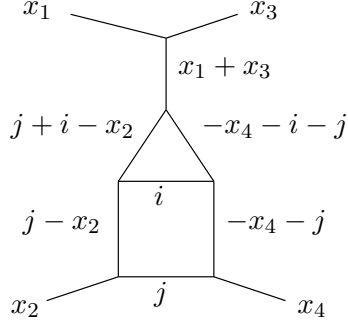

In order to distinguish these labels from the usual weights, let us denote them by $w'(e)$.

The hyperplanes $w'(e) = 0$ subdivide $\mathbb{R}^2$ into different chambers that we call $F$-chambers (see Definition 2.15). In the interior of each $F$-chamber $w'(e)$ has a given sign, and thus every edge inherits an orientation. This also determines the actual weight $w(e) = |w'(e)| = \pm w'(e)$ of the edge. Figure 3 shows the hyperplane arrangement and the corresponding directed graphs with the induced orientations. Since the orientation of the ends and the edge with label $x_1 + x_3$ does not depend on $i$ and $j$, we do not include these edges in the pictures.

Only the bounded $F$-chambers (shaded) correspond to directed **x**-graphs that contribute to the Hurwitz number. The unbounded $F$-chambers correspond to graphs with a directed cycle, hence with multiplicity 0.

For different chambers, the product $\varphi_\Gamma$ differs at most by the sign, since the edge weights $w(e)$ equal plus or minus the edge label $w'(e)$, depending on the side of the hyperplane $w'(e) = 0$ the $F$-chamber is situated. Thus we can define a sign for each $F$-chamber that is determined by the number of edges that are reversed when compared to the reference orientation.

We now develop a formalism that generalizes this discussion.

2.3. **Hyperplane arrangements: formalities.** An **x**-graph $\Gamma$ with a fixed reference orientation can be viewed as a one dimensional cell complex. The differential $d : \mathbb{R}E_\Gamma \to \mathbb{R}V_\Gamma$, sending a directed edge to the difference of its head and tail vertices, gives a short exact sequence:

$$(3) \qquad 0 \to \ker(d) \to \mathbb{R}E_\Gamma \to \operatorname{im}(d) \to 0.$$

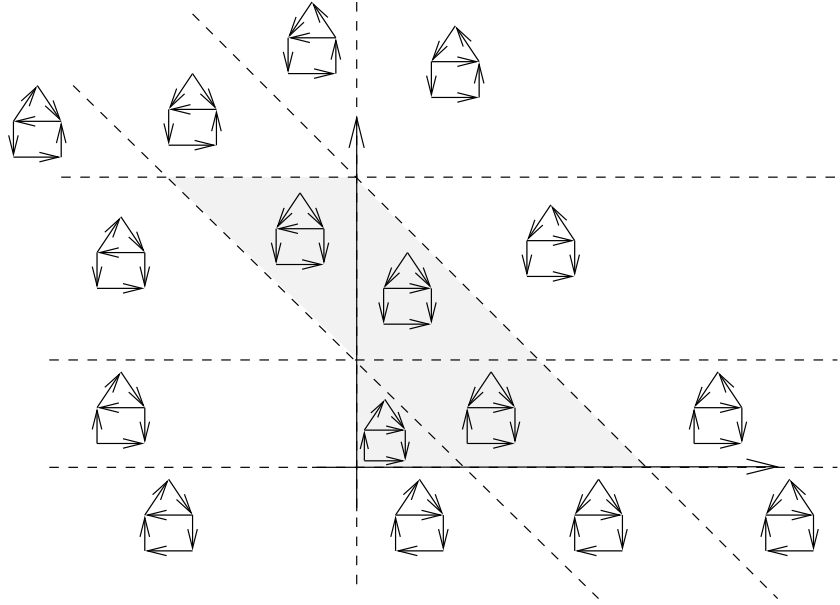

Figure 3. The parameter space for monodromy graphs corresponding to a given **x**-graph.

Decomposing the space of vertices $\mathbb{R}V_\Gamma = \mathbb{R}^n \oplus \mathbb{R}^{i.v.}$ into ends and internal vertices, a vector $(\mathbf{x}, \mathbf{0})$ lies in the image of $d$ when $\sum x_i = 0$.

**Definition 2.8.** The *space of flows* is
$$F_\Gamma(\mathbf{x}) = d^{-1}(\mathbf{x}, \mathbf{0}).$$
Inside it, we have the hyperplane arrangement
$$\mathcal{A}_\Gamma(\mathbf{x})$$
given by the restriction of the coordinate hyperplanes corresponding to the $3g-3+n$ internal edges in $\mathbb{R}E_\Gamma$. The defining polynomial for this hyperplane arrangement is
$$\varphi_\mathcal{A} = \prod_{i=1}^{3g-3+n} e_i,$$
where $e_i$ are the coordinate functions on $\mathbb{R}E_\Gamma$ restricted to $F_\Gamma(\mathbf{x})$.

*Remark* 2.9. Note that $\mathcal{A}_\Gamma(\mathbf{x})$ is a hyperplane arrangement only for generic choices of **x**. In Example 2.5, if **x** satisfies $x_1 + x_3 = 0$, then the space of flows $F_\Gamma(\mathbf{x})$ is contained in the restriction of the hyperplane of the edge from vertex 1 to 2. We still speak of a hyperplane arrangement and consider this a nontransversality of the hyperplane arrangement - a single hyperplane that intersects in codimension 0 rather than 1 as expected in a transverse intersection.

If $x_1 + x_3 \neq 0$, the hyperplane of this edge does not meet the space of flows at all. In general, for an edge which is not part of a cycle, the weight and thus the orientation is determined by $\mathbf{x}$, (same proof as Lemma 6.4 of [CJM]). The corresponding coordinate hyperplane either does not intersect the space of flows, or contains it.

*Remark* 2.10. Note that $F_\Gamma(\mathbf{0}) = H^1(\Gamma, \mathbb{R})$ is a $g$ dimensional vector space; the other $F_\Gamma(\mathbf{x})$ are thus $g$-dimensional affine spaces modeled on $H^1(\Gamma, \mathbb{R})$.

*Remark* 2.11. Our hyperplane arrangement is a variation of a standard construction in algebraic combinatorics. From the short exact sequence (3) there are two natural central hyperplane arrangements, obtained from taking the coordinate hyperplanes in $\mathbb{R}E_\Gamma$, and either taking their image in $\mathrm{im}(d)$ or restricting them to $\ker(d)$. In the literature, these are referred to as the *graphic* and *cographic arrangement*, respectively.

Thus, $\mathcal{A}_\Gamma(\mathbf{0})$ is the cographic arrangement, and $\mathcal{A}_\Gamma(\mathbf{x})$ is a deformation of the cographic arrangement obtained by translating the hyperplanes.

The chambers of $\mathcal{A}_\Gamma(\mathbf{x})$ are indexed by orientations of edges of $\Gamma$, with not all orientations appearing. An orientation occurs if and only if it admits a flow where water is conserved, i.e. if there are no sources or sinks, not even after contracting a cycle.

The chambers we care about are bounded (Section 2.2) and correspond to orientations without cycles (Remark 2.2). These facts are related:

**Lemma 2.12.** *The bounded chambers of $\mathcal{A}_\Gamma(\mathbf{x})$ correspond to orientations of $\Gamma$ with no directed cycles.*

*Proof.* Suppose that a given flow has a directed cycle. Then we may add any fixed positive integer to the weight of each edge in this directed cycle and the balancing conditions are still met and none of the signs of the edges change. Thus, the chamber containing this flow is unbounded.

Now suppose we have an unbounded chamber $A$; then $A$ contains some flow $f$ and some edge $e$ so that the flow along $e$ is greater than $deg = |\mathbf{x_0}|$. Because of the balancing conditions, it is clear that in this case $e$ must be part of a directed cycle in $A$, for water is conserved and only $deg$ enters and leaves the graph. □

**Corollary 2.13.** *Given an **x**-graph $\Gamma$, the bounded chambers of $\mathcal{A}_\Gamma(\mathbf{x})$ are in bijection with directed **x**-graphs projecting to $\Gamma$ after forgetting the orientations of the edges that come from a monodromy graph.*

We conclude this section by establishing some more notation that allows to rephrase equation (2) for $H_g(\mathbf{x})$ in a more convenient way.

**Definition 2.14.** $S_\Gamma(\mathbf{x})$ denotes the contribution to $H_g(\mathbf{x})$ of all monodromy graphs having underlying **x**-graph $\Gamma$.

**Definition 2.15.** For an **x**-graph $\Gamma$, we call $F$-*chambers* the chambers of $\mathcal{A}_\Gamma(\mathbf{x})$ in the flow space $F_\Gamma(\mathbf{x})$. For an $F$-chamber $A$, let $\Gamma_A$ denote the directed **x**-graph $\Gamma$ with the edge directions corresponding to $A$. We use $m(A)$, or $m(\Gamma_A)$, to denote the number of orderings of the vertices of $\Gamma_A$. By Lemma 2.12, $m(A)$ is zero if and only if $A$ is unbounded. We use $\mathrm{Ch}(\mathcal{A}_\Gamma(\mathbf{x}))$ to denote the set of chambers of $\mathcal{A}_\Gamma(\mathbf{x})$.

The sign of $\varphi_\mathcal{A}$ alternates on adjacent $F$-chambers (since we swap the direction of one edge); we use $\mathrm{sign}(A)$ to denote the sign of $\varphi_\mathcal{A}$ on the chamber $A$.

**Definition 2.16.** For integer values of $\mathbf{x}$, the space of flows $F_\Gamma(\mathbf{x})$ has an affine lattice, coming from the integral structure on $\mathbb{Z}E_\Gamma$. We denote this lattice
$$\Lambda = F_\Gamma(\mathbf{x}) \cap \mathbb{Z}E_\Gamma.$$

This notation allows for a convenient interpretation of $S_\Gamma(\mathbf{x})$ in terms of the hyperplane arrangement $\mathcal{A}_\Gamma(\mathbf{x})$. Choices of weights of the edges — i.e. the choice of a flow $f$ on $\Gamma$— correspond to lattice points in $\Lambda$. The product of all the edge weights of a flow $f$ is the absolute value of $\varphi_\mathcal{A}(f)$, which if $f \in A$ is $\mathrm{sign}(A)\varphi_\mathcal{A}(f)$. Thus, we have that

$$(4) \qquad S_\Gamma(\mathbf{x}) = \frac{1}{\mathrm{Aut}(\Gamma)} \sum_{A \in \mathrm{Ch}(\mathcal{A}_\Gamma(\mathbf{x}))} \mathrm{sign}(A) m(A) \sum_{f \in A \cap \Lambda} \varphi_\mathcal{A}(f).$$

## 3. Piecewise Polynomiality

In this section we use Equation (4) to show that double Hurwitz numbers are piecewise polynomial (in Section 3.1), determine the walls (Section 3.2), and show that the polynomials are odd/even (Section 3.3). Section 3.4 contains a discussion of the strong Piecewise polynomiality conjecture of [GJV05].

### 3.1. Polynomials.

**Theorem 3.1** ([GJV05]). *$H_g(\mathbf{x})$ is a piecewise polynomial of degree $4g - 3 + n$.*

*Proof.* The proof is immediate from the following fact: summing a polynomial of degree $d$ over the lattice points in an $g$-dimensional integral polytope of fixed topology is a polynomial of degree $d + g$ in the numbers defining the boundary of the polytope.

This is an analogue of integration of polynomials over a region, and can be seen by iterated applications of Bernoulli's formula for the sum of the first $n$ $k$-th powers. The key point to be careful about is that our vertices are always integers - otherwise one gets quasipolynomials instead of polynomials. Since $\varphi_\mathcal{A}$ is a polynomial of degree $3g - 3 + n$ and $\mathbf{x}$ and $F_\Gamma(\mathbf{x})$ has dimension $g$, the contribution $S_\Gamma(\mathbf{x})$ is locally a polynomial of degree $4g - 3 + n$, and so $H_g(\mathbf{x})$ is as well. $\square$

3.2. **Walls.** From the discussion in the Section 3.1, the functions $S_\Gamma(\mathbf{x})$ are polynomial as long as the topology of the arrangement $\mathcal{A}_\Gamma(\mathbf{x})$ does not change. If we could translate the hyperplanes of $\mathcal{A}_\Gamma(\mathbf{x})$ independently, then generically there would only be transverse intersections, and the topology would change exactly as we passed through nontransverse intersections. In our case, certain nontransversalities occur for every value of $\mathbf{x}$ - however, it is still true that the topology of $\mathcal{A}_\Gamma(\mathbf{x})$ changes when there are additional nontransversalities.

We call the nontransversalities occurring for every $\mathbf{x}$ *good*. The good nontransversalities are easily described: at every interior vertex, if the flow at any two of the adjacent edges is zero, then by the balancing condition the flow at the third edge must also be zero. Thus, for each vertex we have three hyperplanes intersecting in codimension two. The only good nontransversalities are these intersections and their consequences. More explicitly:

**Definition 3.2.** Suppose a set $I$ of $k$ hyperplanes (equivalently, edges in $\Gamma$) in $\mathcal{A}_\Gamma(\mathbf{x})$ intersect in codimension $k - \ell$. We call this intersection *good* if there is a set $L$ of $\ell$ vertices in $\Gamma$ so that $I$ is precisely the set of edges incident to vertices in $L$.

*Remark* 3.3. Recall from Remark 2.9 that we also consider the case where the space of flows is contained in a coordinate hyperplane of an edge a nontransversality. It is not a good nontransversality.

**Definition 3.4.** The *discriminant locus* $\mathcal{D} \subset \mathbb{R}^n$ is the set of values of $\mathbf{x}$ so that for some directed $\mathbf{x}$-graph $\Gamma$ the hyperplane arrangement $\mathcal{A}_\Gamma(\mathbf{x})$ has a nontransverse intersection that is not good. The discriminant is a union of hyperplanes, which we call the *discriminant arrangement*.

The chambers of the discriminant arrangement are the chambers of polynomiality for Hurwitz numbers.

**Definition 3.5.** We call the hyperplanes defining the discriminant arrangement *walls*. The chambers of polynomiality for Hurwitz numbers are called *H-chambers*.

We prove that the walls correspond to the resonance hyperplanes.

**Definition 3.6.** A *simple cut* of a graph $\Gamma$ is a minimal set $C$ of edges that disconnects the ends of $\Gamma$; i.e., there are two ends of $\Gamma$ such that every path between them contains an edge of $C$, and this is true of no proper subset of $C$.

For an **x**-graph $\Gamma$, a flow in $F_\Gamma(\mathbf{x})$ is *disconnected* if for some simple cut $C$ the flow on each edge of $C$ is zero.

*Remark* 3.7. If a flow is disconnected, it follows by the balancing condition that the sum $\sum_{i \in I} x_i$ of weights of ends belonging to a connected component of $\Gamma \setminus C$ is $0$.

**Lemma 3.8.** *The discriminant arrangement $\mathcal{D}$ is given by the set of $\mathbf{x} \in \mathbb{R}^n$ such that for some $\mathbf{x}$-graph $\Gamma$, $F_\Gamma(\mathbf{x})$ admits a disconnected flow.*

*Proof.* Let $f \in F_\Gamma(\mathbf{x})$ be a disconnected flow, and let $C = \{e_1, \ldots, e_k\}$ be a simple cut for $f$. Let $H_1, \ldots, H_k$ be the corresponding hyperplanes; we claim that the $H_i$ intersect nontransversely. The affine space $\bigcap H_i$ corresponds to the space of flows on $\Gamma \setminus C$ with ends **x**. This is an affine space modeled on $H^1(\Gamma \setminus C, \mathbb{R})$. Since we removed $k$ edges from $\Gamma$, the Euler characteristic increased by $k$. But the number of connected components of $\Gamma$ increased by at least one, and so we see that $H^1(\Gamma \setminus C, \mathbb{R})$ has dimension at least $g - k + 1$, and so the $H_i$ intersect nontransversely. If this nontransversality was good, there must be a vertex such that all three adjacent edges belong to $C$. But this contradicts the fact that $C$ is a simple cut, because we do not have to cut all three edges adjacent to a vertex to disconnect.

Now suppose $\mathbf{x} \in \mathcal{D}$. Let $H_1, \ldots, H_k \subset F_\Gamma(\mathbf{x})$ be a maximal collection of hyperplanes having bad non-transverse intersection: $K = \bigcap H_i$ is not contained in any other hyperplane. Since the intersection is nontransverse, by the above reasoning removing the corresponding set of edges $C = \{e_1, \ldots, e_k\}$ disconnects $\Gamma$.
**Claim:** If all the ends of $\Gamma$ lie on the same component of $\Gamma \setminus C$, the other components must consist of a collection of vertices.
This claim implies that $K$ is a good nontransversality, a contradiction. Thus, all **x** on the discriminant arrangement have disconnected flows.

To prove the claim let $e$ be an edge in $\Gamma \setminus C$ not in the component containing all the ends; for any flow $f \in K$, then $f$ must be zero along

$e$, hence $K$ is contained in the corresponding hyperplane, which by maximality must be one of the $H_i$. □

**Theorem 3.9.** *The walls for the discriminant arrangement are given by the resonances, i.e. by equations*
$$\sum_{i \in I} x_i = 0$$
*for any proper subset $I \subset \{1, \ldots, n\}$.*

*Proof.* We have established that for each graph $\Gamma$, the walls of polynomiality of $S_\Gamma(\mathbf{x})$ are the set of $\mathbf{x}$ so that $\Gamma$ admits a disconnected flow. By Remark 3.7, this is a subset of the resonance arrangement. But for any given resonance, it is easy to construct a $\Gamma$ realizing the given resonance: take a graph $\Gamma$ with some edge $e$ so that $\Gamma \setminus e$ has two components, one containing the ends of $I$ and one the other containing the ends of $I^c$. □

3.3. **Parity.**

**Theorem 3.10.** $H_g(\mathbf{x})$ *is either odd or even.*

We first give an elementary proof of this theorem. A more conceptual approach, which plays an important role later, is discussed in Remark 3.11.

*Proof.* Recall that the polynomial $H_g(\mathbf{x})$ is obtained as a sum of polynomials one for each appropriate directed graph. Each of these polynomials is obtained by iterated applications of formulas for sums of $k$-th powers applied to the homogeneous polynomial $\varphi_\mathcal{A}$.

The sum over lattice points for a general even or odd polynomial is typically neither even nor odd; one obtains essentially Bernoulli polynomials, which are neither even nor odd. However, they are nearly so: in the formula for the sum of $k$-th powers, all nonzero terms have the parity of $k+1$ except for the coefficient of $n^k$, which is $1/2$, independently of $k$. This independence means that if $p(x)$ is, say, an odd polynomial, then
$$\sum_{i=1}^n p(i) = q(n) + \frac{1}{2} p(n),$$
where $q$ is an even polynomial. Thus, the error from being even or odd is a constant times the original polynomial on the boundary of the lattice polytope. The point is that $\varphi_\mathcal{A}$, the polynomial we add over lattice points, is the product of the defining equations of the hyperplanes, and so vanishes on the boundary of the polytope. Thus at each step the error is zero, and the resulting polynomials are either odd or even. □

*Remark* 3.11. Conceptually, Theorem 3.10 is a consequence of Ehrhart reciprocity. For $\Delta \subset V$ a $g$-dimensional polytope in a vector space $V$, $\varphi$ a homogenous polynomial of degree $d$ in $t$ and the coordinates of $V$, and $t$ a positive integer, the sum of $\varphi$ over the lattice points in the $t$-dilate of $\Delta$ is polynomial in $t$, which we denote $P_{\Delta,\varphi}(t)$:

$$P_{\Delta,\varphi}(t) = \sum_{x \in \Lambda \cap t\Delta} \varphi(x).$$

Based upon our intuition from integrals, for $t$ a negative integer, we might expect this polynomial to be $(-1)^{g+d}$ times the sum over the lattice points in $-t\Delta$. Ehrhart reciprocity tells us this intuition is not quite correct: the closed polytope $\Delta$ must be replaced by the corresponding open polytope $\Delta°$ - that is, we only sum over the *interior* lattice points of $-t\Delta$:

(5) $$P_{\Delta,\varphi}(-t) = (-1)^{g+d} P_{-\Delta°,\varphi}(t).$$

In our application we are using $\varphi = \varphi_\mathcal{A}$, which is homogenous of degree $3g - 3 + n$. The key point is that $\varphi_\mathcal{A}$, being the defining polynomial of the hyperplane arrangement, vanishes on the boundary of all the polytopes, and so we have

(6) $$P_{\Delta,\varphi_\mathcal{A}}(t) = P_{\Delta°,\varphi_\mathcal{A}}(t).$$

Thus, the correction required by Ehrhart reciprocity vanishes in our case, and Equations 5 and 6 together say our polynomial has parity $4g - 3 + n$.

3.4. **Lower degree bound.** In [GJV05], the Strong piecewise polynomiality conjecture suggests that the polynomial $H_g(\mathbf{x})$ should only have terms of degree in between $4g-3+n$ and $2g-3+n$. Our methods cannot at this time prove this conjecture (although the second author has recently proved this using the infinite wedge, see [Joh]) because it does not hold graph by graph: some graphs contribute monomials of degree lower than $2g - 3 + n$, which remarkably cancel between graphs. Example 3.12 illustrates this situation.

**Example 3.12.** Let $g = 2$ and $n = 2$. The conjecture states that in $H_2(x_1, x_2)$ there are no monomials of degree lower 3.

Figure 3.12 shows all genus 2 directed **x**-graphs with $\mathbf{x} = (x_1, x_2)$. There is always only one choice of vertex ordering. The automorphisms of the second graph contribute a factor of $\frac{1}{4}$, the one of the third graph a factor of $\frac{1}{2}$. The first contributes

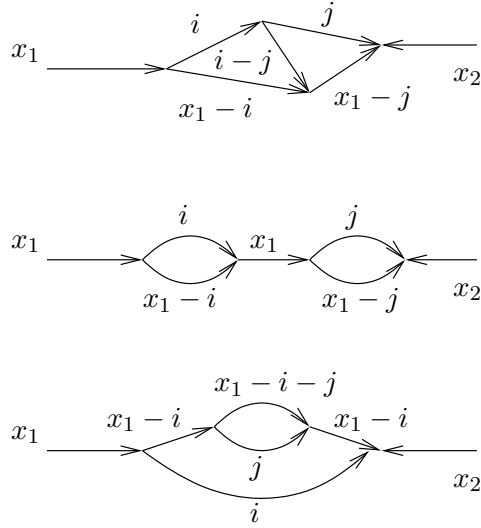

FIGURE 4. All genus 2 directed $(x_1, x_2)$-graphs.

$$f_1 = \sum_{i=0}^{x_1} \sum_{j=0}^{i} i \cdot j \cdot (i-j) \cdot (x_1 - i) \cdot (x_1 - j)$$
$$= \frac{1}{280}x^7 - \frac{1}{60}x^5 + \frac{1}{120}x^3 + \frac{1}{210}x,$$

which has a monomial of degree 1.

The second and third contribute

$$f_2 = \frac{1}{4} \cdot \sum_{i=0}^{x_1} \sum_{j=0}^{x_1} i \cdot j \cdot (x_1 - i) \cdot x_1 \cdot (x_1 - j)$$
$$= \frac{1}{144}x^7 - \frac{1}{72}x^5 + \frac{1}{144}x^3, \text{ resp.}$$
$$f_3 = \frac{1}{2} \cdot \sum_{i=0}^{x_1} \sum_{j=0}^{x_1 - i} i \cdot j \cdot (x_1 - i)^2 \cdot (x_1 - i - j)$$
$$= \frac{1}{504}x^7 - \frac{1}{90}x^5 + \frac{1}{72}x^3 - \frac{1}{210}x.$$

In total, we get

$$H_2(x_1, x_2) = f_1 + f_2 + f_3 = \frac{3x^7 - 10x^5 + 7x^3}{240},$$

and the degree 1-monomials that show up in the contributions from two of the directed **x**-graphs cancel leaving a polynomial which has monomials of degree 3 as lowest degree.

## 4. Wall Crossing: strategy of proof

This section contains an overview of our strategy of proof for Theorem 1.5. It is largely informal and centered on examples. We first illustrate the core idea of our construction in the simplified scenario of a graph admitting only one cut (Section 4.1), then provide an outline of the proof (4.2) and finally illustrate this strategy with a simple but non-trivial example (4.3).

### 4.1. The Simple Cut Example.

**Example 4.1.** Let $\mathbf{x} = (x_1, x_2, x_3, x_4)$ and $r = 6$ (i.e. $g = 2$). Let $I = \{1, 3\}$, and let $\mathfrak{c}_1$ be an $H$-chamber next to the wall $W_I$ satisfying $x_1 + x_3 \leq 0$. $\mathfrak{c}_2$ is the opposite $H$-chamber.

Figure 5 shows an $\mathbf{x}$-graph $\Gamma$ (with reference orientation) for which the hyperplane arrangement $\mathcal{A}_\Gamma(\mathbf{x})$ has a bad nontransversality at the wall $W_I$, and the hyperplane arrangements $\mathcal{A}_\Gamma(\mathbf{x}_1)$ and $\mathcal{A}_\Gamma(\mathbf{x}_2)$ over points $\mathbf{x}_1$ and $\mathbf{x}_2$ on opposite sides of the wall. The nontransversality at the wall consists of the three red hyperplanes meeting in codimension 2. On one side of the wall, these three hyperplanes form a simplex which vanishes when we hit the wall. We call it a *vanishing F-chamber*. A new simplex reappears on the other side of the wall, called an *appearing F-chamber* (see Definition 6.1). The directed $\mathbf{x}$-graph corresponding to the appearing chamber has flows from top to bottom, but none from bottom to top, and so can only be realized when $x_1 + x_3 \geq 0$, i.e. on side "2" of the wall, or in $\mathfrak{c}_2$. This gives a general criterion to see from the graphs whether an $F$-chamber is vanishing/appearing or not (see Lemma 6.6). The 6 neighboring chambers appear on both sides of the wall. In the picture, the directions of the three red edges in each of the bounded $F$-chambers is marked. Also, each bounded $F$-chamber is labeled with a letter, and with its (signed) multiplicity.

We want to understand the contribution of $\Gamma$ to the wall-crossing $P_2(\mathbf{x}_2) - P_1(\mathbf{x}_2)$. To understand the contribution to $P_2$, we sum the polynomial $\varphi_\mathcal{A}$ (weighted with sign and multiplicity) over the lattice points in each of the chambers $A$, $B$, ..., $G$. For the polynomial $P_1$, we have to play the same game with the chambers $B'$, ..., $H'$ on top, however, we evaluate this polynomial now at the point $\mathbf{x}_2$ which is not in $\mathfrak{c}_1$ but in $\mathfrak{c}_2$. Thus, we have to "carry" the chambers $B', \ldots, H'$ over the wall, i.e. we need to interpret the region bounded e.g. by the defining hyperplanes of $B'$ on the other side of the wall in terms of the chambers $A, \ldots, G$.

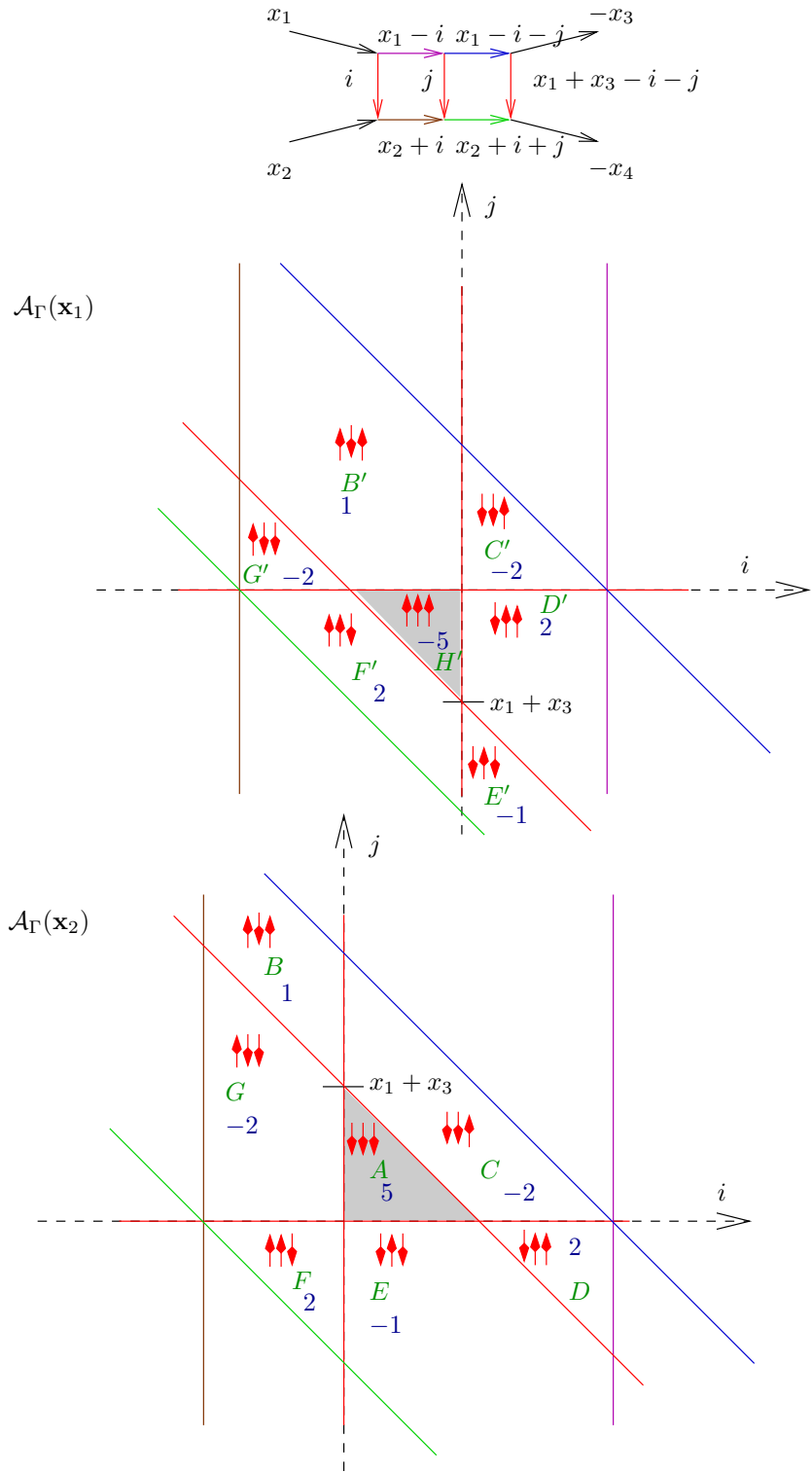

FIGURE 5. An **x**-graph with reference orientations and the hyperplane arrangements for two points on opposite sides of a wall.

We express each of the chambers $B', \ldots, H'$ as a formal signed sum of the chambers $A, \ldots, G$. For example, $C'$ on side 1 is bounded by

$$j \geq 0, x_1 - i - j \geq 0, i \geq 0.$$

The region that is described by these inequalities on side 2 is $A + C$. The hyperplanes bounding simplex $H'$ on side 1 bound $A$, preserving orientation, hence $H' = A$ (recall that we sum the polynomial $\varphi_\mathcal{A}$ over the lattice points of $H'$; we need to switch the summation index twice since $x_1 + x_3 \leq 0$ on side 2, getting a factor of $(-1)^2$).

Altogether, the result of carrying $F$-chambers on side 1 to side 2 can be expressed as follows:

$$H' \to A \qquad B' \to B - A \qquad C' \to C + A \qquad D' \to D - A$$
$$E' \to E + A \qquad F' \to F - A \qquad G' \to G + A.$$

In Section 5 we see that we can view the bounded chambers as a basis of a certain relative homology group of the hyperplane arrangement, and the map just described is the Gauss-Manin connection written in these bases.

Observe that the only chamber on side 1 which contains $B$ in its support when interpreted on side 2 is $B'$, and it contributes positively. Thus, in the difference $P_2(\mathbf{x}_2) - P_1(\mathbf{x}_2)$ the two summands $\sum_B 1 \cdot \varphi_\mathcal{A} - \sum_B 1 \cdot \varphi_\mathcal{A}$ cancel. In fact, all the contributions from chambers which are not appearing chambers cancel, and we only have the contribution from $A$:

$$\sum_A \big(5 - (-5) + 1 - (-2) + 2 - (-1) + 2 - (-2)\big) \varphi_\mathcal{A}$$
(7)
$$= \sum_A 20 \cdot \varphi_\mathcal{A} = \sum_A \binom{6}{3} \cdot \varphi_\mathcal{A}.$$

If we cut the graph $\Gamma$ at the three edges, then the upper part $\Gamma_u$ contributes to the Hurwitz number $H^3(x_1, x_3, -i, -j, -x_1 - x_3 + i + j)$ and the lower part $\Gamma_l$ to $H^3(x_2, x_4, i, j, -x_2 - x_4 - i - j)$. In fact the pair $(\Gamma_u, \Gamma_l)$ appears 6 times in the product of Hurwitz numbers, corresponding to all ways of labelling the three cut edges. Then note that to compute the pair of Hurwitz numbers we must sum over all $i \geq 0$, $j \geq 0$ and $x_1 + x_3 - i - j \geq 0$ (the simplex $A$) the product of internal edges of the two connected components times the connecting edges, hence just the polynomial $\varphi_\mathcal{A}$. Then the contribution to the right hand side of Equation (1) by pair of graphs that glue to $\Gamma$ is $6 \sum_A \binom{6}{3} \cdot \frac{\varphi_\mathcal{A}}{6}$, i.e. Equation (7).

We want to take this a little further, and interpret this equality geometrically. The factor $\binom{6}{3}$ counts the ways to merge two orderings of the vertices of $\Gamma_1$ and $\Gamma_2$ to a total ordering of all vertices. Then re-gluing the cut graphs with the extra data of this merging gives a bijection with the directed, vertex-ordered graphs, contributing to Equation (7).

4.2. **Outline of Proof. Step 1: the left hand side.** For $\mathbf{x} \in \mathfrak{c}_2$, given an $\mathbf{x}$-graph $\Gamma$ and an $F$-chamber $A$, the contribution of $\Gamma$ with edges directed according to the chamber $A$ to $P_2(\mathbf{x})$ is obtained by summing the polynomial $m(A)\varphi_\mathcal{A}$ over $A$. The polynomial $P_1(\mathbf{x})$ is obtained by interpreting in the flow space over $\mathbf{x} \in \mathfrak{c}_2$ the sums and bounds used to compute the Hurwitz numbers in $\mathfrak{c}_1$; thus the graph corresponding to $A$ can contribute to $P_1(\mathbf{x})$: for each chamber $A'$ in $\mathcal{A}(\mathfrak{c}_1)$, consider all hyperplanes that bound $A'$ and trace them in $\mathcal{A}(\mathfrak{c}_2)$: if such hyperplanes bound a set of $F$-chambers including $A$, then assign $A'$ an appropriate sign. Otherwise let $A'$ count 0. Call this coefficient $\langle A', \nabla^*_{\Gamma,12}(A)\rangle$; this choice of notation is explained in Section 5. Then the contribution of $(\Gamma, A)$ to $P_1$ is obtained by summing the polynomial $\varphi_\mathcal{A}$ over $A$ and then multiplying by the number $\left(\sum_{A'} m(A')\langle A', \nabla^*_{\Gamma,12}(A)\rangle\right)$.

**Step 2: the right hand side and the heavy formula.** We again focus on the contribution by the graph $\Gamma$ directed as in the $F$-chamber $A$, and declare a certain subset of edges of $\Gamma$ to be *cuttable*. Ideally we would like to consider all possible ways of cutting $\Gamma$ along cuttable edges into at most three components. We observe that by doing so we recover all graphs contributing to the products of Hurwitz numbers on the RHS of (1) that glue back to $\Gamma$. The polynomial contribution from each cutting of $\Gamma$ is always the same (summing $\varphi_\mathcal{A}$ over $A$), so we would like to show that the signed multiplicity is exactly the one computed on the LHS because of a natural bijection between regluings of the cut graphs (where we allow to reorient the cut edges) and the orientations of $\Gamma$ corresponding to either $A$ or chambers $A'$ such that $\langle A', \nabla^*_{\Gamma,12}(A)\rangle \neq 0$. Alas, cutting the graph into only three connected components doesn't give us enough flexibility to create such a correspondence. Therefore we allow to cut the graphs in all possible legal (see Definition 6.4) ways along cuttable edges, and then organize the inclusion-exclusion process in terms of the number of connected components we have cut the graph into. We therefore wish to prove a wall crossing formula in terms of products of arbitrarily many Hurwitz numbers, that we call the *heavy formula*.

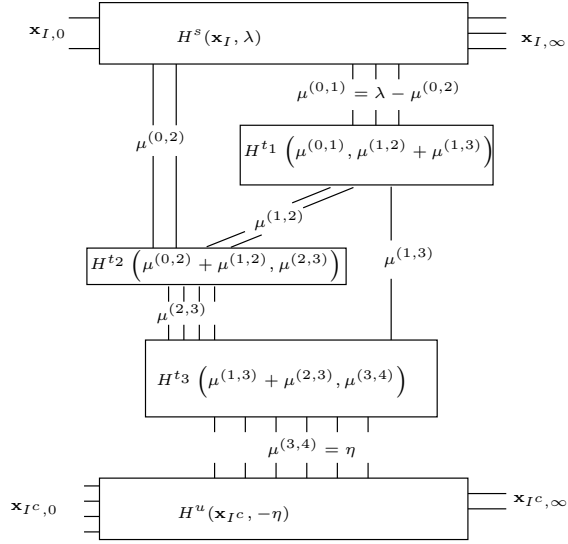

FIGURE 6. The data denoted by $\star$ in the heavy formula, Theorem 4.2

**Theorem 4.2** (Heavy Formula).

$$
(8) \quad WC_I^r(\mathbf{x}) = \sum_{N=0}^{\infty} \sum_{\substack{s+(\sum_{m=1}^N t_m)+u=r \\ |\lambda|=|\eta|=d \\ data\ in\ \star}} \left( (-1)^N \cdot \binom{r}{s, t_1, \ldots, t_N, u} \cdot \right.
$$
$$
\left. \frac{\prod \mu_k^{(i,j)}}{\prod \ell(\mu^{(i,j)})!} \cdot H^s(\mathbf{x}_I, \lambda) \cdot \left( \prod_{m=1}^N H^{t_m}(\star) \right) \cdot H^u(\mathbf{x}_{I^c}, -\eta) \right)
$$

The data denoted by $\star$ is illlustrated in Figure 6: it consists in disconnecting a graph with the right numerical invariants in all possible legal ways, where legal means that the graph obtained by shrinking all connected components to vertices and maintaining the cut edges as edges has no directed cycles. The $\mu_k^{(i,j)}$ denote the partitions of weights of the edges connecting the $i$-th to the $j$-th connected component.

The equivalence of Theorem 4.2 to Theorem 1.5 follows from an inclusion-exclusion argument in Lemma 8.3.

**Step 3: from chambers to cones.** We wish to prove that the multiplicities coming from the Gauss-Manin connection and from the inclusion-exclusion process from the heavy formula match by giving a geometric correspondence between various decorated graphs appearing on the left and right hand side of the formulas. While this is simple "case by case", it is hard to systematize this check, especially because we don't have an efficient expression for the Gauss-Manin connection on the natural basis given by individual chambers. The key observation is that any chamber can be obtained as a linear combination of cones. A cone corresponds to a partial orientation of $\Gamma$, such that cutting along the oriented edges does not disconnect the ends of $\Gamma$; we go further and use this as the definition of a "combinatorial cone" in the space $\mathcal{O}_\Gamma$ of all possible orientations of edges of $\Gamma$. There is a natural lift of a cone to a combinatorial cone, defining a section from the space of $F$-chambers over $\mathfrak{c}_1$ to $\mathcal{O}_\Gamma$. The Gauss-Manin connection acts on combinatorial cones as the identity, and it therefore factors as the composition of the section mentioned above with the natural projection to $F$-chambers over $\mathfrak{c}_2$. We prove Theorem 4.2 by defining a "graph connection" in terms of cuttings and regluings of graphs (see 7.2), that on the one hand acts like the identity on combinatorial cones (and hence agrees with the geometric Gauss-Manin connection), on the other is a natural extension of the inclusion-exclusion in the right hand side of the heavy formula (see 7.3). In this purely combinatorial cutting and gluing process one introduces graphs with sinks or sources, or that lie on the wrong side of the wall. We conclude the proof by checking (Lemmas 7.4, 7.5 and 7.6) that such non-geometric regluings give vanishing contributions, hence this extension recovers the original inclusion-exclusion multiplicity.

4.3. **Following the proof in one example. Step 1: the left hand side.**

**Example 4.3.** Refer to Example 2.7, where we assumed that $0 > x_2 + x_4$. The topology of $\mathcal{A}_\Gamma(\mathbf{x})$ changes if $0 = x_2 + x_4$. Fix the wall $W_{\{2,4\}}$ and let $\mathfrak{c}_1$ and $\mathfrak{c}_2$ be two adjacent $H$-chambers. Assume that in $\mathfrak{c}_1$, we have $0 < x_2 + x_4$, and in $\mathfrak{c}_2$, we have $x_2 + x_4 < 0$. Figure 7 shows the hyperplane arrangements $\mathcal{A}_\Gamma(\mathbf{x}_1)$ and $\mathcal{A}_\Gamma(\mathbf{x}_2)$ for two points $\mathbf{x}_1 \in \mathfrak{c}_1$ and $\mathbf{x}_2 \in \mathfrak{c}_2$. The hyperplanes appear with their defining equations. The bounded $F$-chambers are labeled with letters. Since the edge with weight $x_1 + x_3$ gives the inequality $x_1 + x_3 > 0$, which holds in $\mathfrak{c}_2$ but not on $\mathfrak{c}_1$, every $F$-chamber on the right is an appearing chamber, and every $F$-chamber on the left is vanishing (see Definition 6.1 and Remark 6.2). This is also seen from the corresponding graphs: since over $\mathfrak{c}_2$ the top most interior edge points down, there is a flow from top to bottom (see Lemma 6.6). Figure 8 shows the directed $\mathbf{x}$-graphs corresponding to some of the $F$-chambers.

As in Example 4.1, we pick an appearing $F$-chamber on the right, e.g. $A$, and ask ourselves what $F$-chambers on the left contain it in their support when carried over the wall (this is what we formally define as $\nabla^*_{\Gamma,12}(A)$ in Definition 5.1). To do this, we take chambers on the left, e.g. $E$, and carry them over, i.e. we first determine $\nabla_{\Gamma,12}(E)$. When we carry $E$ over, we get $B$ and keep the orientation, similarly to the vanishing chamber $H'$ in Example 4.1. In the same way, we get $\nabla_{\Gamma,12}(F) = A$. If we interpret the bounds of $G$ in $\mathcal{A}_\Gamma(\mathfrak{c}_2)$, we obtain $A+B+C$ but with reversed orientation. Thus $\nabla_{\Gamma,12}(G) = -A-B-C$. Finally, $H$ becomes $D+B+C$. Thus, $\nabla^*_{\Gamma,12}(A) = F - G$, $\nabla^*_{\Gamma,12}(B) = E - G + H$, $\nabla^*_{\Gamma,12}(C) = -G + H$ and $\nabla^*_{\Gamma,12}(D) = H$.

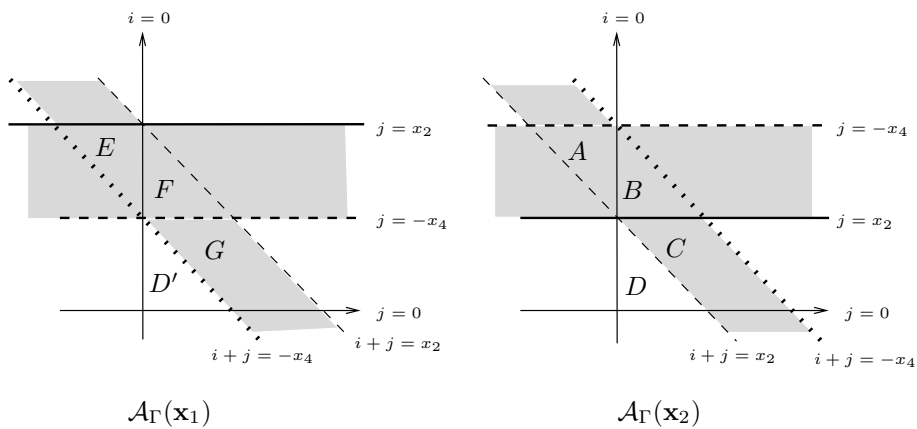

FIGURE 7. The hyperplane arrangements $\mathcal{A}_\Gamma(\mathbf{x}_1)$ and $\mathcal{A}_\Gamma(\mathbf{x}_2)$ for two points $\mathbf{x}_1$ and $\mathbf{x}_2$ on opposite sides of a wall.

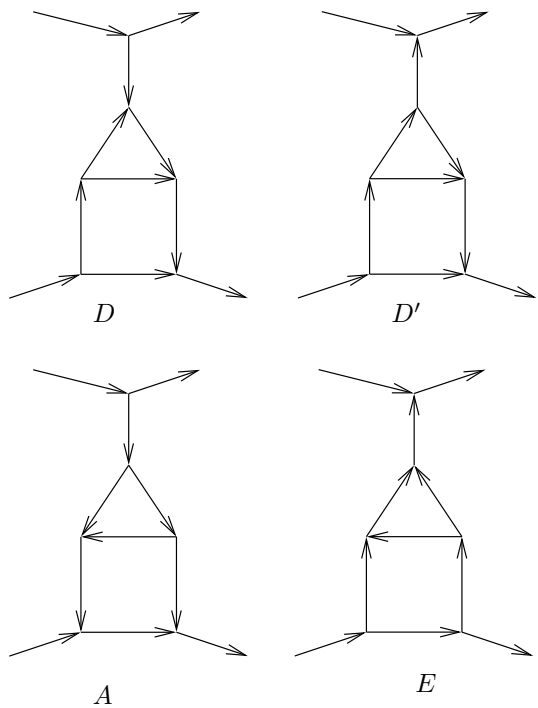

FIGURE 8. The directed $\mathbf{x}$-graphs corresponding to the $F$-chambers $B$, $E$, $F$, $G$ and $H$ of Figure 7.

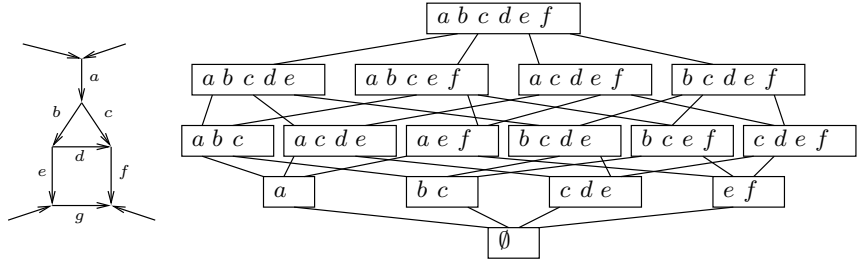

FIGURE 9. The directed **x**-graph $\Gamma_B$ from 4.3 and its poset of $\{1,3\}$-cuts. Letters in the boxes correspond to cut edges.

**Step 2: the right hand side and the heavy formula.** Now we want to establish a bijection between regluings of a cut graph $\Gamma_A$ and graphs in chambers $A'$ with $\langle A', \nabla^*_{\Gamma,12}(A)\rangle \neq 0$.

*Remark* 4.4. The precise statement for this bijection is in Lemma 6.12. Roughly, this lemma states that the number $\langle A', \nabla^*_{\Gamma,12}(A)\rangle$ equals the weighted number of ways to cut the graph $\Gamma_A$ and reglue it to the graph $\Gamma_{A'}$. Each cut is weighted by its rank in a poset of cuts (see 6.4).

**Example 4.5.** We preview the formal definition of a cut with the example of the poset of $\{1,3\}$-cuts of the graph $\Gamma_B$, see Figure 9.

Now we demonstrate the statement of Lemma 6.12. Consider the appearing chamber $B$ in Example 4.3. We determined that $\nabla^*_{\Gamma,12}(B) = E - G + H$. Here are some checks for the weighted number of ways to cut $\Gamma_B$ and reglue it to a given graph:

$\Gamma_E$: to get $\Gamma_E$ from a cut of $\Gamma_B$, we have to turn around the edges $a$, $b$, $c$, $d$, $e$ and $f$. So only cuts that cut all of these edges contribute to $E$. There is only one such cut, the maximal cut. Its rank is 4. We have to turn around all the edges we cut, which is 5 edges. So we get
$$(-1)^4 \cdot (-1)^6 = 1 = \langle E, \nabla^*_{\Gamma,12}(B)\rangle.$$

$\Gamma_G$: to get $\Gamma_G$, we need to turn around $a$, $b$, $c$ and $e$. There are 3 cuts that cut these edges, namely *abcde* and *abcef* (of rank 3), and

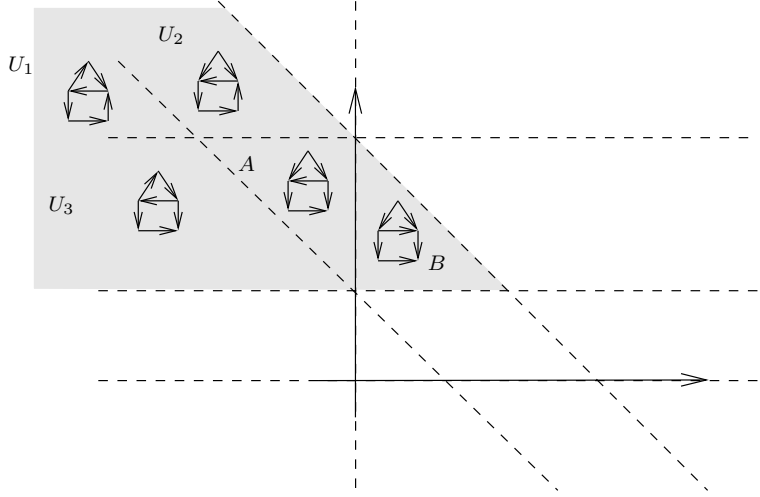

FIGURE 10. A cone

  $abcdef$ of rank 4. Each time, we have 4 edges to turn around. We get:
$$((-1)^3 + (-1)^3 + (-1)^4) \cdot (-1)^4 = -1 = \langle G, \nabla^*_{\Gamma,12}(B) \rangle.$$

$\Gamma_H$: to get $\Gamma_H$, we have to turn around the edges $a$, $b$ and $e$. The cuts that cut these edges are: $abcde$ and $abcef$ of rank 3 and $abcdef$ of rank 4. We have to turn around 3 edges. So we get:
$$((-1)^3 + (-1)^3 + (-1)^4) \cdot (-1)^3 = 1 = \langle H, \nabla^*_{\Gamma,12}(B) \rangle.$$

$\Gamma_F$: to get $\Gamma_F$, we have to turn around $a$, $b$, $c$, $e$ and $f$. The cuts that cut these edges are $abcef$ of rank 3 and $abcdef$ of rank 4. Thus we get
$$((-1)^3 + (-1)^4) \cdot (-1)^4 = 0 = \langle F, \nabla^*_{\Gamma,12}(B) \rangle.$$

**Step 3: from chambers to cones.** To prove our bijection between cutting and regluing of graphs and Gauss-Manin contributions of chambers in Lemma 6.12, we introduce cones.

**Example 4.6.** Figure 10 shows the cone given by:
$$j - x_2 \geq 0 \text{ and } -x_4 - i - j \geq 0$$
in Example 4.3. Combinatorially, the cone requires two edges to have a certain orientation, as depicted in Figure 11. The corresponding combinatorial cone (defined in Section 7) is the set of all graphs for which these two edges have the required orientation, and consists of $2^5 = 32$ oriented graphs. However, 25 of these graphs have a source or sink, as shown in Figure 12. The picture also shows two graphs that don't appear in Figure 10, as they correspond to vanishing chambers. Finally, we have 5 graphs belonging to the cone in Figure 10.

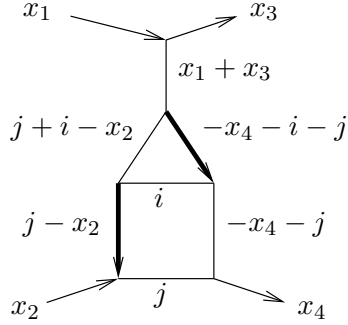

FIGURE 11. The orientations forced by a cone

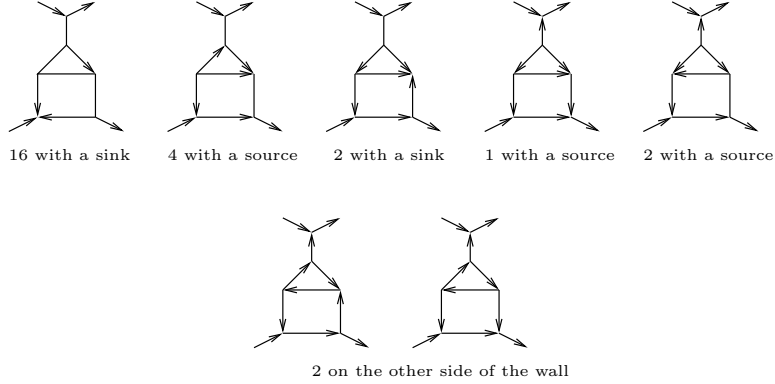

FIGURE 12. A combinatorial cone

Now transport the cone to side one. On this side, the cone consists of only one $F$-chamber, namely the chamber corresponding to the same graph (up to reversing the orientation of the edge $x_1 + x_3$) as chamber $U_1$. We call it $U_1'$. Since the Gauss-Manin connection preserves cones (see Section 7), the sum of all $F$-chambers in the cone $\mathfrak{c}_1$ that map to each $F$-chamber in the cone in $\mathfrak{c}_2$ is one. This is trivially true in this example, since there is only one chamber in the cone on the left, and it maps to each chamber in the cone on the right. We define a graph connection in terms of cutting and regluing of graphs (see Definition 7.2), and show in 7.3 that it acts like the Gauss-Manin connection as the identity on cones. We verify this statement for the example:

- for a chamber inside the cone in $\mathfrak{c}_2$, e.g. $B$: the number of ways to cut $\Gamma_B$ and reglue it so that we stay inside the cone is one.

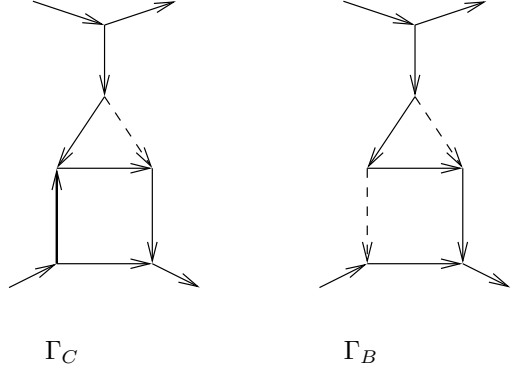

Figure 13. Colored edges for $\Gamma_C$

- for a chamber outside the cone in $\mathfrak{c}_2$, e.g. $C$: the number of ways to cut $\Gamma_C$ and reglue it so that we stay inside the cone in zero.

Remember that the cone fixes the orientation of two edges. In Figure 13, we use different colors for the edges: blue if it is not an edge oriented by the cone, green if it is an edge oriented as prescribed by the cone, red if it points in the wrong direction.

Each nontrivial cut of $\Gamma_B$ e.g. has to cut a blue edge, since the cone edges do not disconnect the ends of the graph. Fix a nontrivial cut, and let $e$ be a blue edge contained in this cut. Now we can pair up regluings in the cone where we keep the orientation of $e$ with regluings where we reverse the orientation of $e$ (and which is also in the cone, since it is a blue edge). The contribution from each pair is 0 since the two graphs differ by the orientation of one edge. Thus the number of ways to nontrivially cut $\Gamma_B$ and stay inside the cone is 0. Since we also have the empty cut, we can altogether cut and reglue exactly once and stay inside the cone. For $\Gamma_C$, essentially the same argument works, only now the empty cut is not a cut which allows us to stay inside the cone, since we have to reverse the red edge.

## 5. Moving hyperplanes and the Gauss-Manin connection

The bounded chambers of the hyperplane arrangement $\mathcal{A}_\Gamma(\mathbf{x})$ can be viewed as a basis for the relative homology group $H_g(F_\Gamma(\mathbf{x}), \mathcal{A}_\Gamma(\mathbf{x}))$, that is:
$$\mathbb{R}[\mathcal{BC}_\Gamma(\mathbf{x})] = H_g(F_\Gamma(\mathbf{x}), \mathcal{A}_\Gamma(\mathbf{x})).$$

Complexifying $F_\Gamma(\mathbf{x})$ and the arrangement $\mathcal{A}_\Gamma^\mathbb{C}(\mathbf{x})$, we still have that $\mathbb{R}[\mathcal{BC}_\Gamma(\mathbf{x})] \cong H_g(F_\Gamma(\mathbf{x}) \otimes \mathbb{C}, \mathcal{A}_\Gamma^\mathbb{C}(\mathbf{x}))$, as the map sending $x + iy$ to $x+(1-t)iy$ is a deformation retract of the pair $(F_\Gamma(\mathbf{x})\otimes\mathbb{C}, \mathcal{A}_\Gamma^\mathbb{C}(\mathbf{x}))$ to the pair $(F_\Gamma(\mathbf{x}), \mathcal{A}_\Gamma(\mathbf{x}))$. Then we can allow $\mathbf{x}$ to take on complex values, giving a complex family of hyperplane arrangements, whose real part is our original family of hyperplane arrangements. The benefit of this maneuver is that now the discriminant $\mathcal{D}$ is a **complex** codimension 1 subvariety, and so its complement $\mathbb{C}^n \setminus \mathcal{D}$ is path connected.

The spaces $\mathbb{R}[\mathcal{BC}_\Gamma(\mathbf{x})]$ form a vector bundle over $\mathbb{C}^n \setminus \mathcal{D}$, which we denote by $\mathcal{BC}$. As a homological bundle, this bundle has a natural flat connection known as the Gauss-Manin connection [Voi02]. In fact, in this case the connection is actually trivial [Var87], and so gives a canonical identification of all the $\mathcal{BC}_\Gamma(\mathbf{x})$. For real $\mathbf{x}$ within one $H$-chamber, this is the obvious identification; for real $\mathbf{x}$ in different $H$-chambers, this is the identification illustrated in Step 1 of the outline of the proof (Section 4.2) and in Examples 4.1 and 4.3.

**Definition 5.1.** Given two $H$-chambers $\mathfrak{c}_1$ and $\mathfrak{c}_2$, and an $\mathbf{x}$-graph $\Gamma$,
$$\nabla_{\Gamma,12} : \mathbb{R}[\mathcal{BC}_\Gamma(\mathbf{x}_1)] \to \mathbb{R}[\mathcal{BC}_\Gamma(\mathbf{x}_2)]$$
denotes the Gauss-Manin connection described above. We give the spaces $\mathbb{R}[\mathcal{BC}_\Gamma(\mathbf{x}_1)]$ and $\mathbb{R}[\mathcal{BC}_\Gamma(\mathbf{x}_2)]$ inner products $\langle \cdot, \cdot \rangle_1, \langle \cdot, \cdot \rangle_2$ by declaring chambers to be an orthonormal basis. We denote by $\nabla_{\Gamma,12}^*$ the adjoint of the Gauss-Manin connection:
$$\langle \nabla_{\Gamma,12}(A), B' \rangle_2 = \langle A, \nabla_{\Gamma,12}^*(B') \rangle_1$$

*Remark* 5.2. Note that the Gauss-Manin connection and its adjoint contain equivalent information, namely:

$\nabla_{\Gamma,12}(A)$**:** the vector in $\mathbb{R}[\mathcal{BC}_\Gamma(\mathbf{x}_2)]$ that corresponds to $F$-chambers bounded by the equations defining $A$ in $\mathfrak{c}_1$.

$\nabla_{\Gamma,12}^*(B')$**:** the vector of $F$-chambers in $\mathbb{R}[\mathcal{BC}_\Gamma(\mathbf{x}_1)]$ whose image via the Gauss-Manin connection contains $B'$.

Note that the chambers in the vectors above appear with a sign corresponding to preserving/reversing orientations.

*Remark* 5.3. For our purposes it is enough to understand $\nabla_{\Gamma,12}$ when $\mathfrak{c}_1$ and $\mathfrak{c}_2$ are two $H$-chambers adjacent across the wall $W(I)$.

For a vector $v$ in a $\mathbb{R}[\mathcal{BC}_\Gamma(\mathbf{y})]$, for some fixed $\mathbf{y}$, the Gauss-Manin connection gives a covariant constant section
$$v : \mathbb{C}^n \setminus \mathcal{D} \to \mathcal{BC}.$$
Given a family of $g$-forms $\omega(\mathbf{x}) \in H^g(F_\Gamma(\mathbf{x}))$ that varies holomorphically in $\mathbf{x}$, then the pointwise pairing
$$\langle v, \omega \rangle = \int_v \omega$$

produces a holomorphic function on $\mathbb{C}^n \setminus \mathcal{D}$.

The continuous analogue of our scenario is the following: the family of forms is given by the family of polynomials $\varphi_\mathcal{A}$, which we want to integrate over a formal sum of bounded $F$-chambers. Performing integration gives naturally a polynomial function in the real $H$-chamber where the topology of the $F$-chambers does not change. When the topology changes, the Gauss-Manin connection tells us that we can keep using the integrating polynomial if we adjust the chambers we're integrating over, and gives the precise prescription for this adjustment.

*Remark* 5.4. Gauss-Manin connections are typically used with continuous structures. Luckily, since our polynomial $\varphi_\mathcal{A}$ vanishes on the boundary of each $F$-chamber, we don't have to deal with the subtleties that arise when dealing with a discrete sum over lattice points (see the discussion in the introduction).

We now interpret the wall crossing in terms of the Gauss-Manin connection.

**From Wall Crossing to Gauss-Manin:** the contribution of the **x**-graph $\Gamma$ to the wall crossing $(:=WC[\Gamma])$ is:

$$WC[\Gamma] = \sum_{A \in \mathcal{BC}_\Gamma(\mathbf{x}_2)} WC[\Gamma, A] \left( \sum_{\Lambda \cap A} \varphi_\mathcal{A} \right),$$

where

$$(9) \quad WC[\Gamma, A] = \text{sign}(A) \left( m(A) - \sum_{B \in \mathcal{BC}_\Gamma(\mathbf{x}_1)} m(B) \langle B, \nabla^*_{\Gamma,12}(A) \rangle \right).$$

## 6. Cuts

Since our hyperplane arrangements have a combinatorial interpretation in terms of graphs, we wish to describe the Gauss-Manin connection directly in terms of the combinatorics of the graphs. In Section 6.3, we state Lemma 6.12, which accomplishes this, then show how the Heavy formula (Theorem 4.2) follows from this lemma. The proof of Lemma 6.12 is postponed until Section 7.

Lemma 6.12 is stated in terms of the poset of cuts, which we introduce in Section 6.1, and identify with the face lattice of a certain cone in Section 6.2.

Throughout this section fix a wall $W_I$ and an **x**-graph $\Gamma$ such that its hyperplane arrangement $\mathcal{A}_\Gamma(\mathbf{x})$ in the flow space $F_\Gamma(\mathbf{x})$ has a bad nontransversality for a point **x** at the wall. Let $\mathfrak{c}_1$ and $\mathfrak{c}_2$ be two $H$-chambers on opposite sides of the wall; by our conventions, $\mathfrak{c}_2$ is the chamber satisfying $\sum_{i \in I} x_i > 0$. Anytime we give a subscript to $\mathbf{x}$, we assume the point lies in the corresponding $H$-chamber.

## 6.1. The poset of cuts.

**Definition 6.1.** Let $A$ be an $F$-chamber of $\mathcal{A}_\Gamma(\mathfrak{c}_2)$. If the inequalities that define $A$ define an empty set in $\mathcal{A}_\Gamma(\mathfrak{c}_1)$, then we call $A$ an *appearing F-chamber*. Analogously, we call $F$-chambers in $\mathcal{A}_\Gamma(\mathbf{x_1})$ that do not exist in $\mathcal{A}_\Gamma(\mathbf{x_2})$ *vanishing F-chambers*.

*Remark* 6.2. If over points of the wall the whole flow space is contained in a coordinate hyperplane ($F_\Gamma(\mathbf{x}) \subseteq \{e_i = 0\}$), then all $F$-chambers in $\mathfrak{c}_2$ are appearing chambers.

**Definition 6.3.** For $\Gamma$ a directed graph, and $E$ a subset of edges of $\Gamma$, we denote by $\Gamma/E$ the directed graph obtained by contracting all edges of $E$. That is, the vertices of $\Gamma/E$ are the connected components of $\Gamma \setminus E^c$, and the edges of $\Gamma/E$ are $E^c$.

Note that $\Gamma/E$ can have cycles and multiple edges.

**Definition 6.4.** For a directed $\mathbf{x}$-graph $\Gamma_A$ and a subset $I \subset \{1, \ldots, n\}$, the set $\text{Cuts}_I(\Gamma_A)$ of *I-cuts of* $\Gamma_A$ consists of subsets $C$ of the edges of $\Gamma_A$ so that either the subset is empty (we say, the *empty cut*) or

(1) $\Gamma_A \setminus C$ is disconnected.
(2) The ends of $\Gamma_A$ lie on precisely two components of $\Gamma_A \setminus C$; one containing all ends in $I$, the other $I^c$.
(3) The directed graph $\Gamma/C^c$ is acyclic and has the (vertex corresponding to the) component containing $I$ as an initial element and that containing $I^c$ as a final element.

The set $\text{Cuts}_I(\Gamma_A)$ is given the structure of a poset by inclusion of cut edges, called poset of $I$-cuts. In Corollary 6.10 we show it is ranked by the number of connected components of the graph with the edges in $C$ removed minus one (so that the empty cut has rank zero). The minimal nonempty elements of the poset $\text{Cuts}_I(\Gamma_A)$ are the simple $I$-cuts of Definition 3.6. The collection of all edges belonging to some $I$-cut is the set of *cuttable edges*.

**Example 6.5.** Remember the graph $\Gamma_B$ from Example 4.3. In Example 4.5, we depicted its poset of $I$-cuts for $I = \{1, 3\}$. In Figure 14, we demonstrate that cutting the edges labeled $b$, $c$, $d$ and $e$ is indeed an $I$-cut (it appears in the poset), while cutting the edges $b$, $c$ and $d$ is not, since the edge $d$ leads to a cycle starting and ending at the bottom connected component.

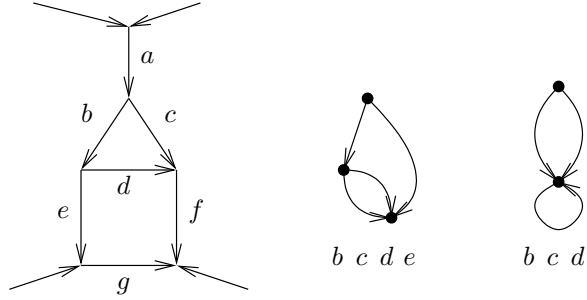

FIGURE 14. The graph $\Gamma_B$ from 4.3 and two cuts.

**Lemma 6.6.** *An F-chamber $A$ in $\mathcal{A}_\Gamma(\mathbf{x_2})$ is an appearing chamber if and only if the directed $\mathbf{x}$-graph $\Gamma_A$ admits an I-cut.*

*Proof.* Assume $\Gamma_A$ admits an $I$-cut $C$. The edges in $C$ are oriented such that they point from the connected component of $\Gamma \setminus C$ containing the ends $I$ to the connected component containing $I^c$. That means there is a flow from the ends $I$ to the ends $I^c$. Such a flow can only exist if $\sum_{i\in I} x_i > \sum_{i\in I^c} x_i$, and this inequality does not hold in $\mathfrak{c}_1$.

Vice versa, if $A$ is an appearing chamber, a subset of the hyperplanes that bound $A$ form a bad nontransversality for a point $\mathbf{y}$ on the wall. Make this nontransversality maximal by adding hyperplanes (not necessarily bounding $A$) containing the non-transverse intersection. Let $C$ be the set of edges corresponding to these hyperplanes.

Lemma 3.8 shows that that $\Gamma_A/C^c$ has a vertex to which all ends in $I$ are contracted (call it $v_I$), a vertex with all ends in $I^c$ and possibly other vertices. The condition that $A$ is an appearing chamber implies that $\Gamma_A/C^c$ has no directed cycles containing both $v_I$ and $v_{I^c}$. We can therefore further contract edges in $\Gamma_A/C^c$ that are contained in some directed cycle (this amounts to shrinking $C$ to a smaller subset $C'$) to obtain an acyclic graph and be assured that we have not identified our two special vertices. The fact that $A \in \mathcal{A}_{\mathfrak{c}_2}$ and $A \notin \mathcal{A}_{\mathfrak{c}_1}$ finally implies that $v_I$ is an initial element and $v_{I^c}$ is terminal. $\square$

6.2. **Geometrization.**

**Definition 6.7.** Consider a directed $\mathbf{x}$-graph $\Gamma_A$ that admits an $I$-cut. Form a new directed graph $\Gamma'_A$ by contracting all vertices above the maximal cut to one vertex, and similarly all vertices below the maximal cut to another vertex. Then $\Gamma'_A$ corresponds to a chamber in the graphical arrangement of $\Gamma'$, defined as the set of $v \in \text{im}(d)$ (the image of $d$ is the subset of $\mathbb{R}V(\Gamma')$ so that the sum of the coordinates are zero) satisfying inequalities $\delta_e(v)(:=$ the difference of the coordinates of $v$ for the tail and source vertices of the edge $e) \geq 0$ for all edges $e \in \Gamma'$.

For $C \in \mathrm{Cuts}_I(\Gamma_A)$ we define a subset $X_C$ of this cone by

$$X_C = \left\{ v \in \mathrm{im}(d) \subset \mathbb{R}V(\Gamma') \,\middle|\, \begin{array}{ll} \delta_e(v) > 0 & e \in C \\ \delta_e(v) = 0 & e \notin C \end{array} \right\}.$$

Furthermore, for a subset $S \subset \mathrm{Cuts}_I(\Gamma_A)$, we use

$$X_S = \bigcup_{C \in S} X_C.$$

Note that with this notation, the entire cone is $X_{\mathrm{Cuts}_I(\Gamma_A)}$, and the vertex of the cone is $X_\emptyset$.

**Definition 6.8.** For any polyhedron $P$, we denote its lattice of faces by $L(P)$.

The following lemma, suggested to us by Federico Ardila, allows us to interpret our inclusion-exclusion type sums as Euler characteristics.

**Lemma 6.9.** *The geometrization mapping $C \mapsto X_C$ induces an isomorphism of posets:*

$$Cuts_I(\Gamma_A) \cong L(X_{Cuts_I(\Gamma_A)}).$$

*Proof.* The elements of both posets are labeled by certain allowable subsets of edges of the graph: for the poset of cuts $\mathrm{Cuts}_I(\Gamma_A)$ these are the cuts, and for the poset of faces $L(X_{\mathrm{Cuts}_I(\Gamma_A)})$ these are the set of all hyperplanes containing a given face. Since in both posets the ordering is given by inclusion, it is enough to show that the allowable subsets of edges for each poset agree.

The poset of cuts $\mathrm{Cuts}_I(\Gamma_A)$ has been carefully discussed already, so we now focus on the poset of faces $L(X_{\mathrm{Cuts}_I(\Gamma_A)})$. Any subset $S$ of edges defines some face: intersect all of the hyperplanes corresponding to edges not in $S$ with the cone $X_{\mathrm{Cuts}_I(\Gamma_A)}$. $S$ is not allowable if the resulting face lies in additional hyperplanes corresponding to edges in $S$; we now understand how this happens.

Suppose then that $v_0 = v_\infty$ is such a hyperplane for the set $S$; that is, $v_0 = v_\infty$ is one of the hyperplanes contained in $S$, but the intersection of the cone and all the hyperplanes not in $S$ is contained in $v_0 = v_\infty$. Since all of the hyperplane equations are of the form $v_i = v_j$ and all of the inequalities defining the cone are of the form $v_i \geq v_j$, the only way to have such an equality forced is to have a cycle of inequalities

$$v_0 \succ v_1 \succ \cdots \succ v_k \succ v_\infty \succ w_1 \succ \cdots \succ w_\ell \succ v_0,$$

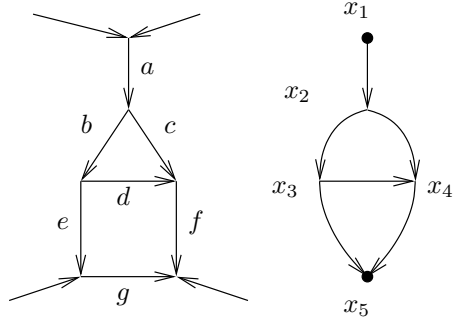

Figure 15. Shrinking everything except the edges of the maximal cut to a vertex.

where $v_i \succ v_j$ means either $v_i = v_j$ is an edge not in $S$ or $v_i = v_j$ is an edge in $S$ and $v_i \geq v_j$ is an inequality defining the cone $X_{\text{Cuts}_I(\Gamma_A)}$. But this is exactly saying there is a cycle in the induced graph of components $\Gamma_A/S^c$; hence the allowable sets of edges in each poset agree.

Similarly, if a set of edges $S$ is not a cut, then there is a directed cycle in the graph $\Gamma_A/S^c$, and the corresponding inequalities force an equality.

From the construction of the graph $\Gamma'_A$ it is clear that a cut must separate the top from the bottom. □

**Corollary 6.10.** *The poset $\text{Cuts}_I(\Gamma_A)$ is ranked, and*

$$rk(C) = \text{the number of components of } \Gamma \setminus C - 1.$$

*Proof.* The poset $L(X_{\text{Cuts}_I(\Gamma_A)})$ is ranked, and the rank is the dimension. From the definition, it is clear that $X_C$ has dimension equal to the number of components of $\Gamma \setminus C$, minus one from restricting to the image of $d$. The statement then follows immediately from Lemma 6.9. □

**Example 6.11.** Return to Example 4.5. Figure 15 shows the graph $\Gamma'_B$ where we shrink everything except the maximal $\{1,3\}$-cut to a vertex. We labeled the vertices of $\Gamma'_B$. These are the coordinates of $\mathbb{R}V_{\Gamma'_B}$. The cone $X_{\text{Cuts}_I(\Gamma_B)}$ is given by the inequalities of the edges, i.e. the edge $a$ gives $v_1 \geq v_2$, $b$ gives $v_2 \geq v_3$ etc. Altogether, we get $v_1 \geq v_2 \geq v_3 \geq v_4 \geq v_5$.

The face lattice of this cone is shown in Figure 16. For a proper face, we turn some of the above inequalities into equalities. Below these equalities, we write down the letters of the edges in $S$ corresponding to the inequalities that we did not turn into equalities. Notice that the lattice is opposite to the poset of cuts in Example 4.5.

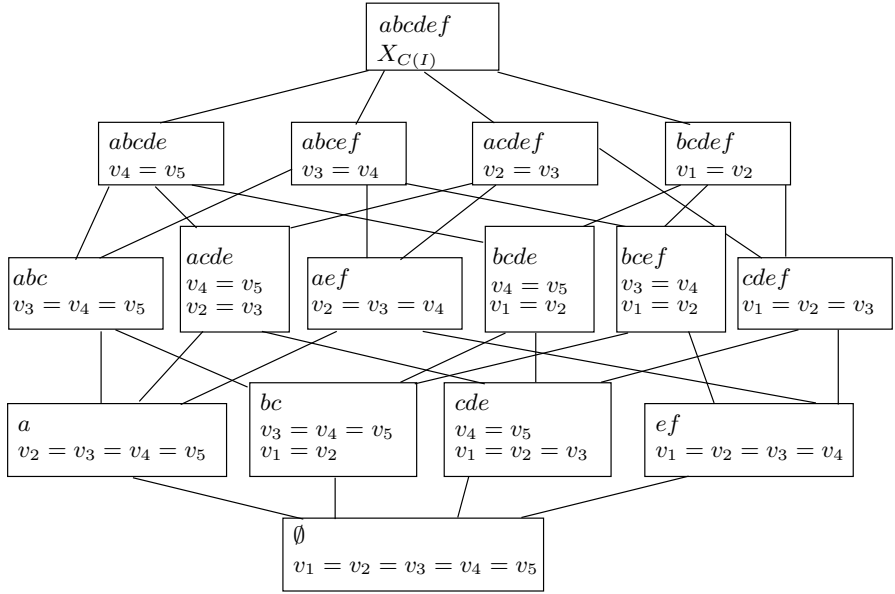

FIGURE 16. The face lattice of $X_{\mathrm{Cuts}_I(\Gamma_B)}$.

6.3. **Combinatorial formula for the Gauss-Manin connection.** The following lemma is the key step connecting the wall crossing expressed through the Gauss-Manin connection with the inclusion-exclusion of products of Hurwitz numbers.

**Lemma 6.12.** *Let $A$ be a (not necessarily bounded) F-chamber in $\mathcal{A}_\Gamma(\mathbf{x}_2)$, and $\Gamma_A$ be the corresponding directed $\mathbf{x}$-graph. For $E$ a subset of the cuttable edges of $\Gamma_A$, consider the graph $\Gamma_{A,E}$ obtained from $\Gamma_A$ by reversing the edges in $E$. If it corresponds to an F-chamber in $\mathcal{A}_\Gamma(\mathbf{x}_1)$, denote this chamber $A_E$. Then:*

$$\nabla^*_{\Gamma,12}(A) = \sum_{C \in Cuts_I(\Gamma_A)} (-1)^{rk(C)} \sum_{E \subset C} (-1)^{|E|} A_E.$$

*If $\Gamma_{A,E}$ has sinks or sources, or if it corresponds to an F-chamber in $\mathcal{A}_\Gamma(\mathbf{x}_2)$ other than $A$, then*

$$\sum_{C \in Cuts_I(\Gamma_A) \ s.t. \ E \subseteq C} (-1)^{rk(C)} = 0.$$

*Remark* 6.13. In case $A$ is not an appearing chamber, then $\mathrm{Cuts}_I(\Gamma_A)$ consists only of the empty cut, and $\nabla^*_{\Gamma,12}(A) = A$, hence $A$ does not contribute to the wall crossing.

*Remark* 6.14. Note that $A_E$ might not be a bounded chamber. Part of the lemma is that these terms contribute zero to the inclusion-exclusion.

For an example see 4.5. We prove Lemma 6.12 in Section 7. Assuming the lemma, we now prove Theorem 4.2.

**From Lemma 6.12 to Theorem 4.2.**

For an $F$-chamber $A$, recall that $m(A)$ counts the number of possible vertex orderings (compatible with the edges directions) for a directed **x**-graph $\Gamma_A$. We can interpret expression (9) as follows: let $\mathfrak{O}(\Gamma)$ denote the set of directed **x**-graphs $\Gamma(d, o)$ with ordered vertices that project to $\Gamma$ when forgetting the extra structure. For $\Gamma(d, o) \in \mathfrak{O}(\Gamma)$, denote by $A_{\Gamma(d)}$ the $F$-chamber identified by the edges directions of $\Gamma(d, o)$. Then:

$$WC[\Gamma, A]\mathrm{sign}(A) = \sum_{\Gamma(d,o)\in\mathfrak{O}(\Gamma)\, s.t.\, A_{\Gamma(d)}\in\mathcal{BC}_\Gamma(\mathbf{x_2})} \langle A_{\Gamma(d)}, A\rangle$$

$$- \sum_{\Gamma(d,o)\in\mathfrak{O}(\Gamma)\, s.t.\, A_{\Gamma(d)}\in\mathcal{BC}_\Gamma(\mathbf{x_1})} \langle A_{\Gamma(d)}, \nabla^*_{\Gamma,12}(A)\rangle,$$

which is nonzero only if $A$ is an appearing chamber.

We now turn to the right hand side of our heavy formula (8), which is an inclusion-exclusion of products of Hurwitz numbers, and it can be computed in terms of tuples of (directed, vertex-ordered) graphs that glue to a graph of genus $g$ with the appropriate number of ends. We isolate the contribution $(:=H[\Gamma])$ by tuples of graphs that reglue to a fixed $\Gamma$. Note that, up a numerical coefficient, the polynomial that we sum over appropriate regions is just $\varphi_\mathcal{A}$, hence:

$$H[\Gamma] = \sum_{A\in\mathcal{BC}_\Gamma(\mathbf{x}_2)} H[\Gamma, A] \left(\sum_{\Lambda\cap A} \varphi_\mathcal{A}\right).$$

Every tuple comes weighted with a multinomial coefficient, corresponding to a merging $\mathfrak{m}$ of the orderings of the connected components to a total order of all vertices. Each such merging produces a graph $\Gamma(d, o) \in \mathfrak{O}(\Gamma)$ by gluing appropriately the ends of the $N + 2$ connected graphs and orienting such edges according to the total order induced by the merging. We can express the contribution to the right hand side of our wall crossing formula by tuples of graphs that reglue to $\Gamma(d, o)$. We denote by $(\gamma_1, \ldots, \gamma_{N+2}, \mathfrak{m} \mid A)$ such a tuple with the extra condition that the "bounds of summation" of the tuple of graphs contain the chamber $A$ (which happens precisely if $A_{\Gamma(d)} = A$ or $A_{\Gamma(d)} = B \in \mathcal{BC}_\Gamma(\mathbf{x_1})$ with $\langle B, \nabla^*_{\Gamma,12}(A)\rangle \neq 0$). Then:

$$H[\Gamma, A] = \sum_{\Gamma(d,o) \in \mathfrak{O}(\Gamma)} \text{sign}(A_{\Gamma(d)}) \sum_{N \geq 0} (-1)^N \sum_{(\gamma_1,\ldots,\gamma_{N+2},\mathfrak{m}|A)} \frac{1}{\prod_{\eta \in \star} \ell(\eta)!},$$

where $\star$ is as in Figure 6. If we forget the order of the ends to glue, each tuple $(\gamma_1, \ldots, \gamma_{N+2}, \mathfrak{m})$ appears in the product of Hurwitz numbers $\prod_{\eta \in \star} \ell(\eta)!$ times (because ends are labeled in our definition of Hurwitz numbers), and the combinatorial factor just cancels such overcounting.

Finally, we need to analyze for what $(\gamma_1, \ldots, \gamma_{N+2}, \mathfrak{m})$ we sum over a chamber $A$. This happens precisely when $(\gamma_1, \ldots, \gamma_{N+2})$ corresponds to a nonempty $I$-cut of $\Gamma_A$. In particular it follows that $A$ is an appearing chamber, since otherwise $\Gamma_A$ would not allow nonempty $I$-cuts (see 6.6). Observing that $N + 1 = \text{rk}(C)$ and that $\text{sign}(A_{\Gamma(d)}) = \text{sign}(A) \cdot (-1)^{|E|}$ we have:

$$H[\Gamma(d,o), A]\text{sign}(A) = - \sum_{C \in Cuts_{\Gamma_A}(I)\setminus\{\emptyset\} \ s.t. \ E \subseteq C} (-1)^{\text{rk}(C)+|E|},$$

where $H[\Gamma(d,o), A]$ denotes the summand corresponding to $\Gamma(d,o)$ in $H[\Gamma, A]$ and $E$ denotes the set of edges for which the orientation differs in $\Gamma_A$ and $\Gamma(d,o)$. Now we observe:

(1) If $A_{\Gamma(d)} = A$ then $E = \emptyset$. It follows from Lemma 6.9 that the ranked sum of all cuts (including the empty cut) is 0, since it is the Euler characteristic of a cone. Since the rank of the empty cut is 0 we have $H[\tilde{\Gamma}, A]\text{sign}(A) = 1$.
(2) If $A_{\Gamma(d)} \neq A$ is in $\in \mathcal{A}(\mathbf{x_2})$ or if $\Gamma(d,o)$ has sinks or sources, then $H[\Gamma(d,o), A]\text{sign}(A) = 0$ by the second part of Lemma 6.12.
(3) If $A_{\Gamma(d)} \in \mathcal{A}(\mathbf{x_1})$, then $E \neq \emptyset$ since $A$ is an appearing chamber and does not exist in $\mathcal{A}(\mathbf{x_1})$. Thus any cut which contains $E$ is not the empty cut and therefore we can include the empty cut in the sum on the right hand side without changing the sum. Then $H[\Gamma(d,o), A]\text{sign}(A) = -\langle A_{\Gamma(d)}, \nabla^*_{12}(A)\rangle$ by the first part of Lemma 6.12.

We have thus shown that:

$$WC[\Gamma, A] = H[\Gamma, A]$$

and hence that Theorem 4.2 follows from Lemma 6.12.

*Remark* 6.15. In our wall crossing formulas Equations (8) and (1), we sum Hurwitz numbers over many different values. For the "middle"

Hurwitz numbers in the expression, these values do not lie in a fixed Hurwitz chamber, so we are not summing a polynomial. However, the left and right Hurwitz numbers (e.g. from Equation (8) the Hurwitz numbers $H^s(\mathbf{x}_I, \lambda)$ and $H^u(\mathbf{x}_{I^c}, -\eta))$ lie in a single chamber, determined the chambers in the original problem:

We cross a wall $\mathbf{x}_I = 0$ and work on the side where $\mathbf{x}_I < 0$. If $H^s(\mathbf{x}_I, \mathbf{y})$ did take values in different walls, we would have $\mathbf{x}_K + \mathbf{y}_J = 0$ for some $\mathbf{x}$ in our chamber and $K \subset I$ and some possible choice of $\mathbf{y}$ and $J$. The sets $K$ and $J$ cannot be empty, as $\mathbf{y}$ is strictly positive and $\mathbf{x}$ lies in the interior of a Hurwitz chamber.

Further since $\mathbf{y}$ positive, this would give $\mathbf{x}_K < 0$ but $\mathbf{x}_K + \mathbf{y} > 0$. But $\mathbf{y} \leq -\mathbf{x}_I$, and so $\mathbf{x}_K - \mathbf{x}_I = -\mathbf{x}_{I \setminus K} > 0$. Thus we have three resonances for $\mathbf{x}$ that we know which side $\mathbf{x}$ was on: $\mathbf{x}_K < 0$, $\mathbf{x}_{I \setminus K} < 0$ and $\mathbf{x}_K + \mathbf{x}_{I \setminus K} = \mathbf{x}_I$. By assumption, our chamber for $\mathbf{x}$ bordered this last wall. But before we can cross this wall, we would have to cross one of the other two walls $\mathbf{x}_K$ or $\mathbf{x}_{I \setminus K}$. Thus we are not adjacent to the wall $\mathbf{x}_I = 0$, which is a contradiction.

## 7. Cones

One of the main insights of [Var87] is that while the Gauss-Manin connection is complicated to write in terms of the bounded $F$-chambers, things simplify when considering unbounded polyhedra. In particular we focus on cones, which are preserved by the Gauss-Manin connection.

The Gauss-Manin connection naturally extends to

$$\nabla_{\Gamma, 12} : \mathbb{R}[\mathrm{Ch}(\mathcal{A}_\Gamma(\mathbf{x}_1))] \to \mathbb{R}[\mathrm{Ch}(\mathcal{A}_\Gamma(\mathbf{x}_2))].$$

By a cone $\mathcal{K} \subseteq F_\Gamma$, we mean a region bounded by hyperplanes $H_i, i \in I$ so that $\cap H_i \neq \emptyset$ for all $\mathbf{x}$. For $\mathfrak{c}$ an $H$-chamber, we use

$$\mathcal{K}(\mathfrak{c}) \in \mathbb{Z}[\mathrm{Ch}(\mathcal{A}_\Gamma(\mathfrak{c}))]$$

to denote the sum of all the chambers in $\mathcal{K}$. Cones are preserved by the Gauss-Manin connection: $\nabla_{\Gamma, 12} \mathcal{K}(\mathfrak{c}_1) = \mathcal{K}(\mathfrak{c}_2)$. Furthermore, cones generate $\mathbb{Z}[\mathrm{Ch}(\mathcal{A}_\Gamma(\mathbf{x}))]$. Since our formula from Lemma 6.12 for the Gauss-Manin connection applies to unbounded chambers, we can prove it by showing that it preserves cones.

A cone $\mathcal{K}$ can be labeled by a partial orientation $\mathcal{P}$ of $\Gamma$: each hyperplane $H_i$ defining $\mathcal{K}$ gives an edge $e_i \in \Gamma$; the orientation of $e_i$ is determined by which side of $H_i$ $\mathcal{K}$ lies on. This labeling is not necessarily unique. For instance, suppose that $e_1, e_2, e_3$ are three edges bordering a vertex $v$, and $H_i$ are the corresponding hyperplanes. Then the $H_i$ split $F_\Gamma(\mathbf{x})$ into six cones, intersecting in the linear space $\cap H_i$.

Each of these cones may be indexed by a partial ordering by only two of the $e_i$, but could also be indexed by an ordering of all three edges.

Any partial orientation $\mathcal{P}$ corresponds to some (perhaps zero) elements of $\mathbb{Z}[\mathrm{Ch}(\mathcal{A}_\Gamma(\mathbf{x}))]$. The first step is to understand which partial orientations $\mathcal{P}$ of $\Gamma$ correspond to cones.

**Lemma 7.1.** *Let $\mathcal{P}$ be a partial orientation of $\Gamma$, with $E$ the set of oriented edges. Then $\mathcal{P}$ defines a cone if and only if all ends of $\Gamma$ lie on the same component of $\Gamma \setminus E$.*

*Proof.* Let $e_i, i \in I$ be the edges of $E$, and $H_i$ the corresponding hyperplanes. For $\mathcal{P}$ to define a cone, the only property we need is that
$$\cap_{i \in I} H_i \neq \emptyset.$$
However, $\cap H_i$ can be understood as the space of flows on $\Gamma \setminus E$. If all the ends of $\Gamma$ lie on one component, then clearly there are flows; however if some component of $\Gamma \setminus E$ contains a proper subset $\emptyset \neq J \neq [n]$ of the ends, then the balancing condition implies $\sum_{j \in J} x_j = 0$, which holds only for $\mathbf{x}$ on the wall $W_J$ and not in general. $\square$

For $\mathcal{P}$ a partial ordering satisfying the conditions of Lemma 7.1, we use $\mathcal{K}_\mathcal{P}$ to denote the corresponding cone.

In fact, our formula for the Gauss-Manin connection does not seem to depend on the geometry at all, only on the directed graphs:

**Definition 7.2.** Let $\mathcal{O}_\Gamma$ denote the set of all orientations of $\Gamma$. Then the *graph connection*
$$\nabla^\mathcal{O}_{\Gamma,12} : \mathbb{R}[\mathcal{O}_\Gamma] \to \mathbb{R}[\mathcal{O}_\Gamma]$$
is defined as follows. Let $A \in \mathcal{O}_\Gamma$ correspond to an oriented graph $\Gamma_A$. For $E$ a subset of the cuttable edges of $\Gamma_A$, denote by $A_E \in \mathcal{O}_\Gamma$ the vector corresponding to having reversed the orientations of the edges in $E$. Then
$$\nabla^{\mathcal{O}*}_{\Gamma,12}(A) := \sum_{C \in Cuts_{\Gamma_A}(I)} (-1)^{\mathrm{rk}(C)} \sum_{E \subset C} (-1)^{|E|} A_E.$$

For any chamber $\mathfrak{c}$, the set $\mathcal{O}_\Gamma$ is the disjoint union of those orientations labeling geometric chambers of $\mathcal{A}_\Gamma(\mathfrak{c})$, and the rest of the orientations, which we denote $\mathcal{NG}_\Gamma(\mathfrak{c})$, for nongeometric.

This induces splittings
$$\mathbb{R}[\mathcal{O}_\Gamma] = \mathbb{R}[\mathrm{Ch}(\mathcal{A}_\Gamma(\mathfrak{c}))] \oplus \mathbb{R}[\mathcal{NG}_\Gamma(\mathfrak{c})].$$

Given a partial orientation $\mathcal{P}$ defining a cone, we define the *combinatorial cone* $\mathcal{K}^\mathcal{O}_\mathcal{P}$ to be the sum of all orientations of $\Gamma$ that agree with $\mathcal{P}$.

For an example of a combinatorial cone, see 4.6.

In the next lemma we show that the graph connection acts as the identity on combinatorial cones.

**Lemma 7.3.** *For any combinatorial cone $\mathcal{K}_\mathcal{P}^\mathcal{O}$, we have*
$$\nabla_{\Gamma,12}^\mathcal{O}(\mathcal{K}_\mathcal{P}^\mathcal{O}) = \mathcal{K}_\mathcal{P}^\mathcal{O}.$$

*Proof.* From the definition of $\nabla_{\Gamma,12}^\mathcal{O}$, we have:

$$\begin{aligned}(10)\quad \langle \nabla_{\Gamma,12}^\mathcal{O}(\mathcal{K}_\mathcal{P}^\mathcal{O}), A \rangle &= \langle \mathcal{K}_\mathcal{P}^\mathcal{O}, \nabla_{\Gamma,12}^{\mathcal{O}*}(A) \rangle \\ &= \sum_{C \in Cuts_{\Gamma_A}(I)} (-1)^{\mathrm{rk}(C)} \sum_{E \subset C} (-1)^{|E|} \langle \mathcal{K}_\mathcal{P}^\mathcal{O}, A_E \rangle.\end{aligned}$$

Consider a non-trivial cut $C$. As the oriented edges of $\mathcal{P}$ do not disconnect the ends of $\Gamma$, but the edges in $C$ do, there must be an edge $e \in C \cap \mathcal{P}^c$.

It follows that
$$\sum_{E \subset C} (-1)^{|E|} \langle \mathcal{K}_\mathcal{P}^\mathcal{O}, A_E \rangle = 0,$$
as the set of all $E \subset C$ can be partitioned into those which contain $e$ and those which don't, with an obvious bijection between the two sets. As $e \notin \mathcal{P}$, adding or subtracting $e$ only changes the sign $(-1)^{|E|}$, and not $\langle \mathcal{K}_\mathcal{P}^\mathcal{O}, A_E \rangle$, and hence the terms for $E$ containing $e$ cancel with those for $E$ not containing $e$.

The only remaining contribution is from the empty cut, and hence equation (10) simplifies to
$$\langle \nabla_{\Gamma,12}^\mathcal{O}(\mathcal{K}_\mathcal{P}^\mathcal{O}), A \rangle = \langle \mathcal{K}_\mathcal{P}^\mathcal{O}, A \rangle,$$
which proves the lemma. $\square$

For an example of the statement, see 4.6. Lemma 7.3 is very close in form to Lemma 6.12, but the adjoint to the graph connection applied to an $F$-chamber $A \in \mathcal{BC}_\Gamma(\mathbf{x}_2)$ a priori contains contributions from (chambers corresponding to) graphs with sources and sinks, and from graphs that do not occur in $\mathcal{BC}_\Gamma(\mathbf{x}_1)$. The next three lemmas prove Theorem 6.12 by showing that these nongeometric chambers do not actually appear in $\nabla_{\Gamma,12}^{\mathcal{O}*}(A)$.

**Lemma 7.4.** *Suppose $A \in \mathcal{BC}_\Gamma(\mathbf{x}_2)$ and $S \in \mathcal{O}_\Gamma$ corresponds to a graph with a source or sink. Then $\langle \nabla_{\Gamma,12}^\mathcal{O}(S), A \rangle = 0$.*

*Proof.* From the definition,
$$\langle \nabla_{\Gamma,12}^\mathcal{O}(S), A \rangle = \sum_{C \in \mathcal{C}(I)} (-1)^{\mathrm{rk}(C)} \sum_{E \subset C} (-1)^{|E|} \langle S, A_E \rangle.$$

Since $S$ is a single chamber, $\langle S, A_E \rangle$ vanishes unless $E = \Delta$, where $\Delta$ is the set of edges where the orientations of $S$ and $A$ differ. Hence,

$$\langle \nabla^{\mathcal{O}}_{\Gamma,12}(S), A \rangle = \sum_{\substack{C \in C_{\Gamma_A}(I) \\ \Delta \subset C}} (-1)^{\mathrm{rk}(C)}$$

We interpret this sum as the Euler characteristic of an appropriate complex. Consider the dual complex $X^{\vee}_{\mathrm{Cuts}_I(\Gamma_A)}$ and the natural identification of its lattice of faces with $\mathrm{Cuts}_I(\Gamma_A)$. Denote by $X_\Delta$ the subcomplex corresponding to cuts containing $\Delta$. Each minimal cut $C_m \supset \Delta$ corresponds to a maximal cell $X_{C_m} \in X_\Delta$.

We now prove that

$$\bigcap_{\substack{C_m \supset \Delta \\ C_m \text{ minimal}}} \overline{X}_{C_m}$$

is positive dimensional. A source or sink of $S$ must occur at a vertex $v$ interior to the set of cuttable edges of $A$, for flipping all cuttable edges of $A$ gives a graph with no sources or sinks. Since the local picture in $A$ at $v$ begins with no sources or sinks, the orientation of (at least) one the edges incident to $v$ is preserved in the sink/source. Such edge $e$ is not part of any minimal cut that produces $S$ from $A$, and hence the ray in $X^{\vee}_{C_{\Gamma_A}(I)}$ corresponding to cutting all cuttable edges except $e$ is contained in $\bigcap \overline{X}_{C_m}$.

Since $X_\Delta$ is a union of a finite number of cones all of whose (multiple) intersections consist of positive dimensional cones, the Euler characteristic is $\chi(X_\Delta) = 0$. □

**Lemma 7.5.** *If $A, A' \in \mathcal{BC}_\Gamma(\mathbf{x}_2)$ are two distinct appearing chambers, then*

$$\langle A', \nabla^{\mathcal{O}*}_{\Gamma,12}(A) \rangle = 0.$$

*Proof.* As in the previous lemma, if we let $\Delta$ denote the set of edges with different orientations in $\Gamma_A$ and $\Gamma_{A'}$:

$$\langle \nabla^{\mathcal{O}}_{\Gamma,12}(A'), A \rangle = \sum_{\substack{C \in C_{\Gamma_A}(I) \\ \Delta \subseteq C}} (-1)^{\mathrm{rk}(C)}$$

Since by 6.9 the weighted sum of all cuts in $C_{\Gamma_A}(I)$ is 0, we prove the equivalent fact that the weighted sum over all cuts that do not cut all of $\Delta$ is 0. We now work with the face complex $X_{C_{\Gamma_A}(I)}$, and denote $X^c_\Delta$ the subcomplex of cuts that do not cut all of $\Delta$. The maximal cells of $X^c_\Delta$ correspond to maximal cuts $C_M$. Again, we conclude the proof by

arguing that
$$\bigcap_{\substack{C_M \not\supseteq \Delta \\ C_M \text{ maximal}}} \overline{X}_{C_M}$$
is positive dimensional. This is because $\Delta$ is disjoint from the set of cuttable edges of $\Gamma_{A'}$. Any nontrivial cut of $\Gamma_{A'}$ is therefore a cut of $\Gamma_A$ corresponding to a positive dimensional face in $\bigcap \overline{X}_{C_M}$. □

**Lemma 7.6.** *If $A \in \mathcal{BC}_\Gamma(\mathbf{x}_2)$ and $A' \notin \mathcal{BC}_\Gamma(\mathbf{x}_1) \cup \mathcal{BC}_\Gamma(\mathbf{x}_2)$ is a bounded chamber, then*
$$\langle A', \nabla^{\mathcal{O}_*}_{\Gamma,12}(A) \rangle = 0.$$

*Proof.* This lemma argues that a graph cannot cross over a wall different than $W_I$ by just flipping cuttable edges for the $I$-wall. Let $\Gamma_A$ be the graph corresponding to chamber $A$. Consider a wall $W_J$ given by a subset $J \neq I$, and assume that in $\mathfrak{c}_2$ the inequality $\sum_{k \in J} x_k > 0$ holds. Restricting to variables in $J \cap I$ or $J \cap I^c$, at least one of the two sums must be still strictly positive. Assume that $\sum_{k \in I \cap J} x_k > 0$.

If one cuts all $I$-cuttable edges, $\Gamma_A$ has two connected components that contain all the ends, $\Gamma_A^I$ containing all the $I$ ends and $\Gamma_A^{I^c}$ all the $I^c$ ends. Since the weight of the $I$-cuttable edges is arbitrarily small when approaching $W_I$, the flow that enters $\Gamma_A^I$ must leave again through the ends $I$. Thus there must be a directed path from a positive end in $I \cap J$ to a negative end in $I \cap J^c$. But this prevents any regluing of the cut graph to identify an appearing chamber on the side of $W_J$ corresponding to $\sum_{k \in J} x_k < 0$.
□

**Example 7.7.** As an example for the statement of Lemma 7.5, consider the directed **x**-graphs $\Gamma_A$, $\Gamma_{A'}$ in Figure 17, corresponding to appearing chambers. We want to see that the weighted number of ways to cut $\Gamma_A$ and reglue it to $\Gamma_{A'}$ is zero. Any cut that allows to reglue $\Gamma_A$ to $\Gamma_{A'}$ must cut the two dotted edges. We have to show that $(\sum_C (-1)^{\text{rk}(C)}) = 0$, where the sum goes over all cuts that cut at least the two dotted edges. Alternatively, since the weighted sum of all cuts is 0 by 6.9, we can show that the sum of all cuts that do not cut both dotted edges is 0.

$\Gamma_A$ admits 6 simple cuts that are also shown in Figure 17. The maximal cuts which do not cut both wrong edges are 123 and 346. Since the poset of cuts equals the face lattice of $X_{\text{Cuts}_I(\Gamma_A)}$, we can describe the cuts which do not cut both wrong edges as a subset $X_\Delta^c$ of this cone. We want to see that $X_\Delta^c$ is contractible, which implies that its Euler characteristic — which equals the weighted sum of cuts that don't cut both wrong edges — is 0. The figure also shows a schematic picture of $X_\Delta^c$, intersected with the sphere (so that rays become points).

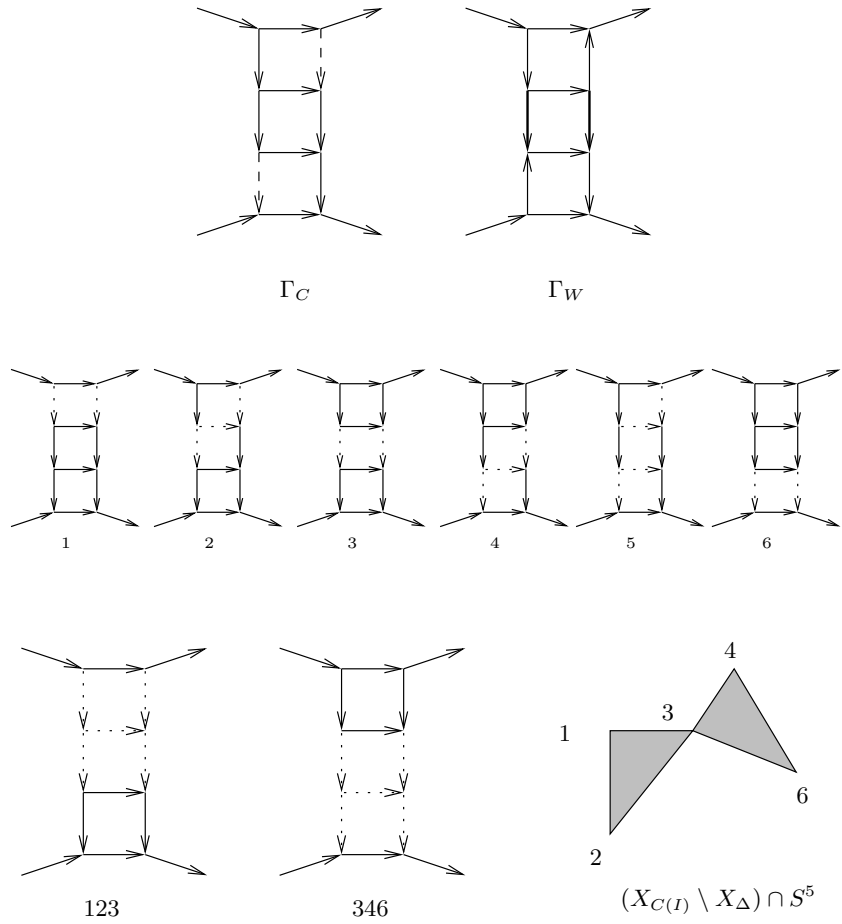

Figure 17. A chamber on the same side of the wall.

We can see that the maximal cells intersect in the ray 3, which is the only cut of $\Gamma_{A'}$. Thus the set is contractible.

## 8. Heavy to Light: Proof of Theorem 1.5

To prove Theorem 1.5 we observe that we are applying the Gauss-Manin connection to a particular vector where each chamber $A$ has coefficient $m(A)\text{sign}(A)$. With these weights, the expression of the Gauss-Manin connection simplifies.

**Definition 8.1.** A cut $C$ where every edge is adjacent to either the top or the bottom component is called a *thin cut*. We refer to the union of connected components that are not the top or the bottom component as the *middle components*.

To each cut $C$, there is a naturally associated thin cut $t(C)$, obtained by keeping the edges in $C$ that border either of the components containing the ends, and forgetting all edges of $C$ that are between two interior components. For each thin cut $T$, the set of cuts $C$ with $t(C) = T$ forms a sub-poset $P(T)$ of the poset of all cuts.

Theorem 1.5 naturally appears as a sum over all thin cuts. To prove it, we match the sum over $P(T)$ of the terms in Theorem 4.2 as an inclusion-exclusion formula in Lemma 8.2 that evaluates to the contribution of cut $T$ in Theorem 1.5. Denote by $\mathfrak{o}(C)$ the number of vertex orderings on $\Gamma \setminus C$.

**Lemma 8.2.** *For a thin cut $T$ with $t$ vertices in the middle components, we have:*

$$(11) \qquad (-1)^t \mathfrak{o}(T) = \sum_{C \in P(T)} (-1)^{rk(C)} \mathfrak{o}(C).$$

To prove Lemma 8.2, we use the following lemma:

**Lemma 8.3.** *Given a thin cut $T$ of a graph $\Gamma$, denote by $\gamma(T)$ the induced directed subgraph on the vertices in the middle components. Using the notation of 6.7, we have*

$$P(T) = L(X_{\gamma(T)}).$$

*Proof.* The proof is a slight modification of that of Lemma 6.9. The same arguments show that those subsets of edges of $\gamma(T)$ that avoid creating directed cycles are exactly those that correspond to faces of $X_{\gamma(T)}$. The only difference is that since every set of edges in $P(T)$ already separates the top from the bottom, we use a different graph, which does not have vertices corresponding to the top and the bottom components. □

*Remark* 8.4. Although every thin cut $T$ is also a cut, the signs with which the corresponding terms appear in Theorems 1.5 and 4.2 need not agree - in Theorem 1.5 the sign of the term corresponding to $T$ is given by the number of components, while in 4.2 it is given by the number of vertices. However, $T$ is the minimal element of $P(T)$, and the maximal element of $P(T)$ corresponds to the cut where each component consists of one vertex, and so the sign of the maximal element of $P(T)$ in the heavy formula does agree with the sign of $T$ in the thin cut formula (see Equation (11)).

Fix a thin cut $T$, and consider the graphical arrangement for the graph $\gamma(T)$; the hyperplanes correspond to edges, and the chambers correspond to orientations of the edges that have no cycles, though they may have sources and sinks. The cone $X_{\gamma(T)}$ is one chamber of this arrangement.

To each chamber $\mathcal{C}$ we associate the number $\mathfrak{o}(\mathcal{C})$ of total orderings of the graph $\Gamma \setminus T$ that are compatible with the directed graph of the top and bottom components of $\Gamma \setminus T$, and with the directed edges of the middle components given by $\mathcal{C}$ (in particular $\mathfrak{o}(X_{\gamma(T)}) = \mathfrak{o}(T)$). Observe that the geometrically opposite chamber $X_{\gamma(T)}^{op}$ corresponds to reversing all the edges of the middle components, and so $\mathfrak{o}(X_{\gamma(T)}) = \mathfrak{o}(X_{\gamma(T)}^{op})$.

For a given cut $C$ in $P(T)$, denote by $\mathcal{K}_C$ the affine tangent cone to $X_{\gamma(T)}$ along the face $X_C$: this is the cone containing $X_{\gamma(T)}$ defined by all facets of $X_{\gamma(T)}$ incident to $X_C$. The cone $\mathcal{K}_C$ contains the chambers corresponding to all possible orientations of the cut edges.

$$\mathfrak{o}(C) = \sum_{\mathcal{C} \in \mathcal{K}_C} \mathfrak{o}(\mathcal{C}).$$

From this discussion, we see that Lemma 8.2 follows from the following identity on cones:

**Lemma 8.5.**
$$(-1)^t X_{\gamma(T)}^{op} = \sum_{C \in P(T)} (-1)^{rk(C)} \mathcal{K}_C$$

*Proof.* Proposition 3.1 of [BV97] gives two inclusion-exclusion formulas for a bounded polyhedron $P$ in terms of cones based at the faces $f \in \mathcal{F}$ of $P$. These formulas are essentially equivalent to our lemma; we must make some adjustments to deal with the fact that we are using unbounded cones instead of bounded polyhedrons.

Before we can state the formulas of [BV97], we must introduce some notation. The inward pointing cones $C_P^+(f)$ are defined by the same hyperplanes that define $f$, and contain $P$, while the outward pointing cones $C_P^-(f)$ is the opposite cone based at $f$; this means that $C_P^+(P)$ and $C_P^-(P)$ are both the whole space. Rather than working with characteristic functions as in [BV97], we take formal sums of chambers. Then we have:

$$(12) \qquad P = \sum_{f \in \mathcal{F}} (-1)^{\dim f} C_P^+(f) \text{ and}$$

$$(13) \qquad (-1)^{\dim P} P = \sum_{f \in \mathcal{F}} (-1)^{\dim f} C_P^-(f).$$

We can prove 8.5 by applying each of (12) and (13) in turn.

To apply these formulas to our situation, take a hyperplane $H$ transverse to $X_{\gamma(T)}$, and let $P = H \cap X_{\gamma(T)}$. Then, the faces of $P$ are in bijection with the **nonempty** cuts in $P(T)$. The cones here are exactly the nonzero inward facing tangent cones; however, because of the shift in the dimension the rank of the cut is the dimension of the face plus one. So, by Equation (12) the sum over the nonempty cuts is giving a contribution of $-X_{\gamma(T)}$, which exactly cancels the contribution from the empty cut.

Let $H_0$ be the hyperplane parallel to $H$ that passes through the origin, and $H_0^+$ the half-space bounded by $H_0^+$ and containing $X_{\gamma(T)}$. Then from the above, we have that $\sum_{C \in P(T)} (-1)^{\text{rk}(C)} \mathcal{K}_C \cap H_0^+ = 0$.

Now translate the hyperplane $H$ past the origin, keeping its normal direction fixed: the resulting $P$ shrinks to a point, and then become a transverse slice of $H \cap X_{\gamma(T)}^{\text{op}}$. The cones in the sum are now the outward pointing tangent cones, and so we essentially have the sum in Equation (13), but again with the signs shifted. This gives us a contribution of $(-1)^{\dim P + 1} X_{\gamma(T)}^{\text{op}} = (-1)^t X_{\gamma(T)}^{\text{op}}$. We have $\sum_{C \in P(T)} (-1)^{\text{rk}(C)} \mathcal{K}_C \cap H_0^- = (-1)^t X_{\gamma(T)}^{\text{op}}$ (where $H_0^-$ denotes the other half-space), and the claim follows.

□

**Example 8.6.** In this example, we demonstrate that the light formula (see 1.5) cannot be proved with a bijection between cut and reglued graphs and chambers contributing to the Gauss-Manin connection as in Lemma 6.12 and Remark 4.4. Consider the graph $\Gamma$ from Figure 18 and the wall $W_I$ with $I = \{1, 2\}$. The wall crossing contribution of this graph equals twice the sum of $\varphi_\Gamma$ over the lattice points in $A$, since $-\langle B, \nabla^*_{\Gamma, 12}(A) \rangle = 1$. Lemma 6.12 tells us that we can interpret this 1 as the minus the number of ways to cut $\Gamma$ and reglue it to $\Gamma_B$, for which all interior edges are reversed. Figure 19 shows all $I$-cuts of $\Gamma$ with the sign corresponding to minus their rank in the poset. We can see that the only cut which allows to reglue to $\Gamma_B$ is the last one, since this is the only cut which cuts all four edges that we need to reverse. So minus the number of ways to cut and reglue to $\Gamma_B$ equals $1 = -\langle B, \nabla^*_{\Gamma, 12}(A) \rangle$. For our light formula, we don't allow to cut edges which are not connected to the top or the bottom part. All the light cuts of $\Gamma$ are depicted in Figure 20. These cuts are weighted with a sign corresponding to the number of vertices of the middle components. We can see that there is no bijection as for the heavy cuts, since there is no way to cut the graph in a light way and still cut the four edges we need to reverse. If we take vertex orderings into account however, we can still recover the $2 = 1 - (-1)$ with which chamber $A$ contributes to the wall crossing: Count each light cut with the binomial factor corresponding to the number of vertex orderings of the disconnected graph. This contribution is

$$\binom{4}{1} + \binom{4}{3} + \binom{4}{2} - \binom{4}{1,1,2} - \binom{4}{2,1,1} + \binom{4}{1,2,1} = 2.$$

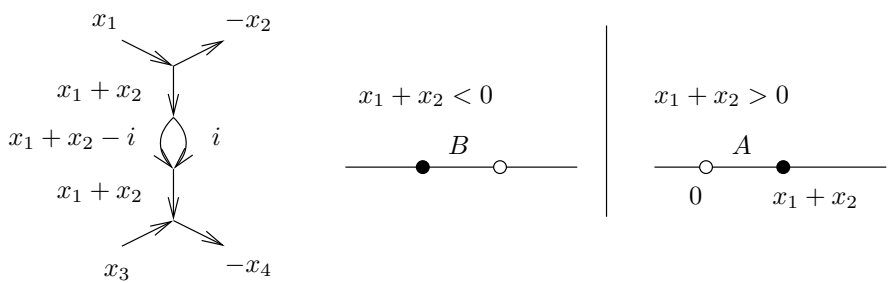

FIGURE 18. An **x**-graph with reference orientation and the two hyperplane arrangements $\mathcal{A}_\Gamma(\mathbf{x}_1)$ and $\mathcal{A}_\Gamma(\mathbf{x}_2)$ left and right of the wall $x_1 + x_2 = 0$.

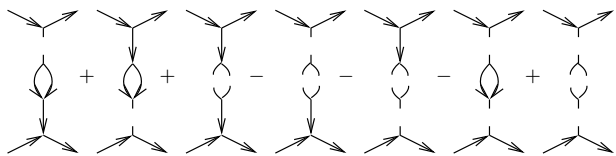

FIGURE 19. All $I$-cuts of $\Gamma$.

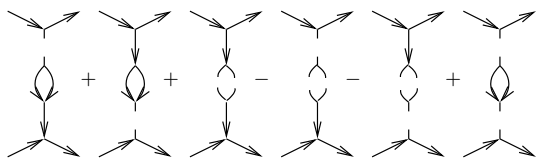

FIGURE 20. All thin cuts of $\Gamma$.

Renzo Cavalieri, Colorado State University, Department of Mathematics, Weber Building, Fort Collins, CO 80523-1874, USA
  *E-mail address*: renzo@math.colostate.edu

Paul Johnson, Department of Mathematics, Imperial College London, 180 Queen's Gate, London SW7 2AZ, UK
  *E-mail address*: paul.johnson@imperial.ac.uk

Hannah Markwig, CRC "Higher Order Structures in Mathematics", Georg August Universität Göttingen, Bunsenstr. 3-5, 37073 Göttingen, Germany
  *E-mail address*: hannah@uni-math.gwdg.de